\documentclass[a4,11t]{article}
\usepackage[latin1]{inputenc}
\usepackage{multicol}
\usepackage{psfrag}
\usepackage{fullpage}
\usepackage[dvips]{graphicx}
\usepackage{color}
\usepackage{epsfig}
\usepackage{deleq}
\usepackage{float}
\usepackage{cite}
\usepackage{amsmath,amssymb,amsfonts}
\usepackage{ntheorem}
\usepackage{cases}
\numberwithin{equation}{section}
\newcommand{\qed}{\hfill $\Box$ \medskip}
\newcommand{\pr}{\noindent {\sc \bf{Proof}:\ }}
\newcommand{\bel}[1]{\begin{equation}\label{#1}}
\newcommand{\be}{\begin{equation}}
\newcommand{\ee}{\end{equation}}
\newcommand{\bea}{\begin{eqnarray*}}
\newcommand{\eea}{\end{eqnarray*}}
\newcommand{\bean}{\begin{eqnarray}}
\newcommand{\eean}{\end{eqnarray}}
\newtheorem{theorem}{Theorem}[section]
\newtheorem{lemma}{Lemma}[section]
\newtheorem{proposition}{Proposition}[section]

\newtheorem{remark}{Remark}[section]
\newtheorem{definition}{Definition}[section]
\newtheorem{cor}{Corollary}[section]

\theoremstyle{plain}
\newcommand{\N}{\mathbb{N}}

\newcommand{\R}{\mathbb{R}}

\newcommand{\ts}{T_*}

\title{ Global Dynamics for a relativistic charged and
colliding plasma in presence of a massive scalar field on the
Robertson-Walker spacetime }
\author{  Marcelin
 Kenmogne Noumo   \\ University of Yaound\'e 1, Faculty of Science\\ Department of
Mathematics, P.O. Box: 812 Yaound\'e Cameroon \\ e-mail: marconoumo@yahoo.fr  \\[2ex]
         Norbert Noutchegueme
                      \\University of Yaound\'e 1, Faculty of Science\\ Department of
Mathematics, P.O. Box: 812 Yaound\'e  Cameroon\\ e-mail: nnoutch@yahoo.fr\\[2ex]
 Roger Tagne Wafo\thanks{Corresponding author}   \\ University of Douala, Faculty of Science\\ Department of
Mathematics and Computer Science, P.O. Box: 24157 Douala  Cameroon \\ e-mail: rtagnewafo@yahoo.com  \\[2ex]
}

\begin{document}
 \maketitle
\begin{abstract}
We consider the coupled Einstein-Maxwell-Boltzmann system with
cosmological constant in presence of a massive scalar field. The
background metric is that of Friedman-Lemaître-Robertson-Walker
space time in the spatially homogeneous case where the unknown
functions only depend on time and not on the space variables
$(x^i)$, $i=1,2,3$. By combining the energy estimates method with
that of characteristics we derive under suitable conditions on the
chock kernel (see (\ref{eq2.22})), a local (in time) solution of the
coupled system. Further, under the hypotheses that the data are
small in some appropriate norms and that the cosmological constant
satisfies $\Lambda
> -4\pi m^2\Phi_0^2$,
 we derive a unique global (in time) solution
(Theorem \ref{MTT}).

\vspace{3mm}
\noindent{\textbf{Keywords}}: Einstein-Maxwell-Boltzmann system;
massive scalar field; Sobolev spaces; local existence; global
existence, continuity argument.
\end{abstract}

\section{Introduction}

\indent The profound knowledge of the universe phenomena is a
scientific preoccupation nowadays. A local and a global modelization
of the universe are then required. The General Relativity built in
1916 by Albert Einstein is essential to understand, explain and
predict the universe phenomena at the macroscopic level. In this
paper, we study the global dynamics of a relativistic plasma in the
Friedman-Lemaître-Robertson-Walker space time. We will only call it
a Robertson-Walker space time. We study the Einstein equations,
which are the basic equations of the General Relativity and which
describe the gravitational forces coupled to the Maxwell equation,
which describe the electromagnetic forces and to the Boltzmann
equation, which is one of the basic equation of the kinetic
relativistic theory, describing the dynamics of the massive and
charged particles, by determining their distribution function $f$,
which is a positive scalar function of the position and the momentum
of the particles. We suppose that we are in presence of a massive
scalar field which can be used to measure the gravitational waves,
which can propagate through the space at the speed of the light,
even in presence of material bodies, analogously to electromagnetic
waves. See \cite{NDM} and references therein for more details on
this subject.
 Now our motivation for considering the Einstein equations with cosmological
constant $\Lambda$, is due to the fact that astrophysical
observations based on luminosity via redshift plots of some far away
objects such as Supernova-Ia, have made evident the fact that the
expansion of the universe is accelerating. A classical mathematical
tool to model this phenomenon is to include the cosmological
constant $\Lambda$ in the Einstein equations. For more details, see
\cite{HL1, NS, NT, NG, Nal, SBT}. Also notice that the scalar fields
are considered to be a mechanism producing accelerated models, not
only in ''inflation'', which is a variant of the Big-Bang theory
including now a very short period of very high acceleration, but
also in the primordial universe.

This work falls in the general framework of relativistic kinetic
theory where the coupled Einstein-Vlasov system and its
generalization, the coupled Einstein-Boltzmann system  are studied.
In \cite{DBYCB}, Daniel Bancel  and Yvonne Choquet-Bruhat, 1973,
proved an  existence and uniqueness theorem of the solution of the
Cauchy problem for the coupled Einstein-Maxwell-Boltzman, the
background metric being a small perturbation of the Minkowski
metric. In \cite{NNEMT}, 2009, the authors proved that the coupled
Einstein-Maxwell-Vlasov system with cosmological constant has a
unique global (in time) solution in  Bianchi I-VIII space-times. In
\cite{HL1, HL2, HLEN}, Ho Lee and Ernesto Nungesser studied the
asymptotic behavior of solutions of the coupled Einstein-Vlasov and
Einstein-Boltzmann systems in Robertson-Walker and Bianchi I
spacetimes. In \cite{NNDD, NT, NNRA}, N. Noutchegueme and his
collaborators considered the Einstein-Maxwell-Boltzmann system on
Bianchi I spacetimes and the Einstein-Boltzmann on Robertson-Walker
spacetimes and proved existence and uniqueness of weak global
solutions  in weighted $L^1$ space.

The paper consider the coupled Einstein-Maxwell-Boltzmann system
with cosmological constant in presence of a massive scalar field.
The background metric is that of Friedman-Lemaître-Robertson-Walker
space time in the homogeneous case where the unknown functions only
depend on time and not on the space variables $(x^i)$, $i=1,2,3$. By
combining the energy estimates method with that of characteristics
we derive under suitable conditions on the chock kernel (see
(\ref{eq2.22})), a local (in time) solution of the coupled system.
Further, under the hypotheses that the data are small in some
appropriate weighted Sobolev norms and that the cosmological
constant satisfies
$\Lambda > -\sigma^2$,
 where $\sigma>0$ is a constant depending only
on the potential of the scalar field, we prove that the coupled
Einstein-Maxwell-Boltzmann system with massive scalar field and
cosmological constant has a unique global (in time) solution
(Theorem \ref{MTT}). This result extends existing results on the
relativistic Boltzmann equation (see for example \cite{NT, HL2,
HLEN, NNAN} or \cite{NNDD} and \cite{NNRA} in a different setting)
for the following:
\begin{itemize}
  \item We study the evolution of matter distribution, the electromagnetic
  field, the scalar field and that of spacetime simultaneously see \cite{HLEN}.
  \item The case of positive cosmological constant is covered by our
  analysis but a range of negative cosmological constants is allowed.
  We should point out the fact that in the absence of matter field,
  by following step by step the proof of Proposition 4.2 page 302 of
  \cite{NT} one can prove that there can not exist global solution in the case
  $\Lambda < 0$. Note that negative $\Lambda$ appears naturally in
   the anti-de Sitter/conformal field theory correspondence and in some
   theoretical models as supergravity, string theory,  etc. See
  \cite{CCN1}.
  \item The approach. Our idea is to combine the geometric techniques of
  energy estimate and the method of characteristics for single
  first order partial differential equations to obtain some
  boundedness properties of the iterative sequence. From these estimates, we derive the
  local weak solution as the limit in some Sobolev spaces of the
  iterative scheme. The main difficulty when using this approach is to
derive from this the classical ($C^1-$) solution of the reduced
system. This is done in two steps using the interpolation inequality
in real order Sobolev spaces and embedding Sobolev inequalities (see
proof of Theorem \ref{MT} on page \pageref{MT}).
   This approach has to be compared with that of
  \cite{NNDD,NT,NNRA} where it is not clear how to apply
  directly the Cauchy-Lipschitz Theorem to the system of ordinary
  differential equations obtained by the method of characteristics.
\end{itemize}

 The paper is organized as follows:
In the second Section, we present in details the equations we are
interested with; namely, the coupled Einstein-Maxwell-Boltzmann
system with massive scalar field and cosmological constant. We
specify the assumptions we impose on the shock kernel of the
collision operator $Q$. Further, the compatibility equations are
given and the equation for the potential $\Phi$ of the scalar field
with positive mass $m$ is derived. In Section 3, we recall a local
existence Theorem for the Boltzmann equation. This gives the
opportunity to introduce the functional spaces we are using and to
recall some Moser-type substitution inequalities for the collision
operator $Q$ which will be use later. In the fourth Section, we
first introduce some new coordinates for which the Boltzmann
equation has a convenient form and secondly
 we
introduce a new set of unknown functions for which the
Einstein-Maxwell-Boltzmann system with massive scalar field become
an equivalent system of first order differential equations; see
Equations (\ref{eq5.12})-(\ref{eq5.18}) and finally the iterative
scheme is set up.
In Section 5, we derive a local existence theorem (Theorem \ref{MT})
to the Einstein-Maxwell-Boltzmann system with massive scalar field
as the limit of a suitable sequence in our functional spaces. In the
last section (Section 6) global existence of solutions to the
Einstein-Maxwell-Boltzmann system with massive scalar field is
established under the hypotheses that the norm of the Cauchy data is
small enough and the Cosmological constant $\Lambda$ is such that
$\Lambda>-\sigma^2$.
The paper ends with an appendix where we derive an energy estimate
for the solution of an hyperbolic first order PDE: the $H^k_d-$norm
of the unknown is estimated in terms of the $H^k_d-$norm  of the
Cauchy data and an integral of the $H^k_d-$norm of the source term.
This is applied to the Boltzmann equation written in some
appropriate coordinates.

\emph{\textbf{ Acknowledgments}}. RTW has the pleasure to thank
Piotr Chru{\'s}ciel for helpful discussions concerning the
continuity on pages 25 and 26. MKN and RTW would like to thank the
Henri Poincaré Institute where part of this work has been carried
out.

\section{The Equations}
Here we give the equations under consideration, then we specify the
simplification hypotheses under which we are studying the coupled
system and finally we study the compatibility of the equations. We
end the section with the study of the conservation conditions from
where we derive the equation for the scalar field $\Phi$.
\subsection{Simplification hypotheses and the detailed equations}
 We are in the flat
Robertson-Walker spacetime $(\R^4,g)$, whose metric $g$ of signature
$(-,+,+,+)$ can be written in the canonic coordinates $(x^{\alpha})$
of $\R^4$ in which $x^0=t$ is the time and $(x^i)$, $i=1,2,3$ the
space variables as:

\begin{equation}\label{eq2.1}
    g=-dt^2+a^2(t)\big[(dx^1)^2+(dx^2)^2+(dx^3)^2\big]
\end{equation}
where $a>0$ is an unknown function called the expansion factor. We
study an homogeneous case, which means that the unknown functions do
not depend on the space variables $(x^i)$, $i=1,2,3$.
We adopt the Einstein summation convention
$A_{\alpha}B^{\alpha}=\sum\limits_{\alpha}A_{\alpha}B^{\alpha}$, the
Greek indices $\alpha$, $\beta$, $\cdots$ vary from 0 to 3, and the
Latin indices $i$, $j$, $k$, $\cdots$ from 1 to 3. The phenomenon is
governed by the following system :
\begin{numcases}
\strut R_{\alpha\beta}-\frac{1}{2}Rg_{\alpha\beta}+\Lambda
g_{\alpha\beta}=8\pi(T_{\alpha\beta}+
\tau_{\alpha\beta}+K_{\alpha\beta}+H_{\alpha\beta}) \label{eq2.2}\\
\nabla_{\alpha}F^{\alpha\beta}=J^{\beta}\label{eq2.3}\\
\nabla_{\alpha}F_{\beta\gamma}+\nabla_{\beta}F_{\gamma\alpha}+\nabla_{\gamma}F_{\alpha\beta}=0 \label{eq2.4}\\
\mathcal{L}_Xf=Q(f,f)\label{eq2.5}
\end{numcases}
where :
\begin{enumerate}
    \item [$\bullet$] $(\ref{eq2.2})$ are the Einstein equations,
    the basic equations of the General Relativity which describe the
    gravitational forces with $\Lambda$ the cosmological constant;
    $R_{\alpha\beta}$ is the Ricci tensor and
    $R=g^{\alpha\beta}R_{\alpha\beta}$ is the Riemann scalar
    curvature. $T_{\alpha\beta}$, $\tau_{\alpha\beta}$,
    $K_{\alpha\beta}$ and $H_{\alpha\beta}$ are the components of
    the stress Mass-Energy tensor, source of the gravitational field
    $g$, with :

\begin{numcases}
\strut
T_{\alpha\beta}=\int_{\mathbb{R}^3}\frac{p_{\alpha}p_{\beta}f(t,\overline{p})}{p^0}a^3dp^1dp^2dp^3\label{eq2.6}\\
\tau_{\alpha\beta}=-\frac{1}{4}g_{\alpha\beta}F^{\lambda\mu}F_{\lambda\mu}+F_{\alpha\lambda}F_{\beta}^{\lambda}
\label{eq2.7}
\\
K_{\alpha\beta}=-\theta_{\alpha\beta}\label{eq2.11}\\
H_{\alpha\beta}=\nabla_{\alpha}\Phi\nabla_{\beta}\Phi-\frac{1}{2}(\nabla^{\lambda}\Phi\nabla_{\lambda}\Phi
+m^2\Phi^2)g_{\alpha\beta}\label{eq2.12}
\end{numcases}
where :
\begin{enumerate}
    \item [$-$] $T_{\alpha\beta}$ is generated by the distribution
    function $f$ of the charged, massive and colliding particles,
    which is a positive scalar function of the time $t=x^0$and the
    momentum $\,p=(p^{\alpha})=(p^0,\overline{p})=(p^0,p^i)$:
    $$
f : T\R^4\simeq\R^4\times\R^4\longrightarrow\R^+,
(x^{\alpha},p^{\alpha})\longmapsto f(x^{\alpha},p^{\alpha})\in\R^+
    $$
    and $a^3=|g|^{\frac{1}{2}}=|\det g|^{\frac{1}{2}}$. We suppose
    that the massive particles move on the mass hyperbolo\"{\i}d
    $(\mathcal{H}) : g(p,p)=-1$, and that they are ejected towards the
    future in the time oriented manifold $\R^4$.\\
    From $(\ref{eq2.1})$, we then deduce that :
     \begin{equation}\label{eq2.8}
p^0=\sqrt{1+a^2\big((p^1)^2+(p^2)^2+(p^3)^2\big)} \;.
    \end{equation}
    Equation $(\ref{eq2.8})$ shows that we can always express $p^0$ in terms of
    $(p^i)$.
    $f$ is solution to the Boltzmann equation $(\ref{eq2.5})$ we
    present later.
    \item [$-$] $\tau_{\alpha\beta}$ is the Maxwell tensor
    associated to the electromagnetic field
    $F=(F_{\alpha\beta})=(F^{0i},F_{ij})$ in which $F^{0i}$ stands
    for the electric part and $F_{ij}$ for the magnetic part. $F$ is
    an antisymmetric and closed $2-$form, solution of the Maxwell
    equations $(\ref{eq2.3})-(\ref{eq2.4})$ and describes the
    electromagnetic forces.
\item[$-$] In Equation (\ref{eq2.11}), $\theta_{\alpha\beta}$ is a symmetric
$2-$tensor called the pseudo-tensor of pressure. The general form of
$\theta_{\alpha\beta}$ is due to A. Lichnerowicz \cite{AL} pages 36
and 37. The cases of pure matter $(\theta_{\alpha\beta}=0)$ and
perfect fluid ($\theta_{\alpha\beta}=p_rg_{\alpha\beta}$, where
$p_r$ is a scalar function representing the pressure) are particular
cases. We make on $\theta_{\alpha\beta}$ the assumptions :
\begin{equation}
\label{eq2.13} \nabla_{\alpha}\theta^{\alpha\beta}= -\rho^2u^{\beta}
\qquad  \mbox{and}\qquad
g^{ij}\theta_{ij}=0
\end{equation}
where in $(\ref{eq2.13}),$ $\rho>0$ is a constant and $(u^{\beta})$
the four velocity of the particles.
\item[$-$] $H_{\alpha \beta}$ is the stress-energy tensor defined by the scalar field $\Phi$ with the positive
mass $m>0$ which is an unknown function of $t$.
\end{enumerate}
    \item [$\bullet$] In $(\ref{eq2.3})$ which stands for the
    first group of the Maxwell equations, $J^{\beta}$ is the
    Maxwell current generated by the charged particles through the
    formula :
    \begin{equation}\label{eq2.9}
J^{\beta}=\int_{\mathbb{R}^3}\frac{p^{\beta}f(t,\overline{p})}{p^0}a^3dp^1dp^2dp^3-eu^{\beta}\;;
\end{equation}
where $e\geq 0$ is an unknown function which stands for the
elementary electric density. We consider that the particles are
comoving, which means $(u^{\beta})=(1,0,0,0)$. In the homogeneous
case we always have $\nabla_{\alpha}F^{\alpha 0}=0$, where $\nabla$
is the Levi-Civita connection of $(\R^4,g)$. So, equation
$(\ref{eq2.3})$ implies $J^0=0$. This determines $e$ to be :
 \begin{equation}\label{eq2.10}
e(t)=a^3(t)\int_{\mathbb{R}^3}f(t,\overline{p})d\overline{p} \;.
\end{equation}
\item [$\bullet$] Equation $(\ref{eq2.4})$ is the second group of the Maxwell
equations and is equivalent to $dF=0$ since $F$ is a closed
$2-$form.
    \item [$\bullet$]  (\ref{eq2.5})  is the Boltzmann equation for
    $f$ we now introduce. In that equation
    $\mathcal{L}_Xf$ is the Lie derivative of $f$ with respect to
    the  vector field $X=X(F)=(p^{\alpha}, \mathcal{P}^{\alpha})$ where :
\begin{equation}\label{eq2.15}
    \mathcal{P}^{\alpha}=-\Gamma^{\alpha}_{\lambda\mu}p^{\lambda}p^{\mu}+eF^{\alpha}_{\lambda}p^{\lambda}
    \;;
\end{equation}
where in this equation, the  $\Gamma_{\;\lambda\mu}^{\alpha}$'s are
the Christoffel symbols of the metric $g$.
The trajectories of the charged particles are curves in $T\R^4 :
s\longmapsto \big(x^{\alpha}(s),p^{\alpha}(s)\big)$ solutions of the
differential system:
\begin{equation}\label{repl}
\frac{dx^{\alpha}}{ds}=p^{\alpha}\qquad \mbox{and}\qquad
\frac{dp^{\alpha}}{ds}=\mathcal{P}^{\alpha} \;;
\end{equation}
where $\mathcal{P}^{\alpha}$ is defined by $(\ref{eq2.15})$.
Equations $(\ref{repl})$ show that
$X=(p^{\alpha},\mathcal{P}^{\alpha})$ is tangent to the trajectories
of the charged particles. Since $f=f(t,p^i)$,  the Boltzmann
equation $(\ref{eq2.5})$ writes :
\begin{equation}\label{eq2.18}
p^0\frac{\partial f}{\partial t}+\mathcal{P}^i\frac{\partial
f}{\partial{p^i}}=Q(f,f) \;.
\end{equation}
We now introduce the collision operator $Q$ which appears in
$(\ref{eq2.18})$.
$Q$ is the binary and elastic operator introduced by A. Lichnerowicz
 and Tchernikov in 1940, and according to which only two
particles enter in collision at a point $(t,x^i)$, without
destructing each other, and if $(p,q)$ are their momenta before the
collision and $(p',q')$ their momenta after the collision, the sums
$p+q$ and $p'+q'$ are preserved. Let $g$ and $h$ be two functions
defined on $\R^3$. We have :
\begin{equation}\label{eq2.19}
    Q(g,h)=Q^+(g,h)-Q^-(g,h)
\end{equation}
with:
\begin{numcases}
\strut
       Q^+(g,
h)=\int_{\mathbb{R}^{3}}\frac{a^{3}(t)}{q^{0}}d\overline{q}\int_{\mathbf{S}^2}g(\overline{p}')
h(\overline{q}')B\big(a,\overline{p}, \overline{q},
\overline{p}',\overline{q}'\big)d\omega \label{eq2.20}\\
       Q^-(g,
h)=\int_{\mathbb{R}^{3}}\frac{a^{3}(t)}{q^{0}}d\overline{q}\int_{\mathbf{S}^2}g(\overline{p})
h(\overline{q})B\big(a,\overline{p}, \overline{q},
\overline{p}',\overline{q}'\big)d\omega\label{eq2.21}  \;.
      \end{numcases}
We now present the different elements of formulae $(\ref{eq2.20})$
and $(\ref{eq2.21})$ point by point, specifying the assumptions we
adopt.
\begin{enumerate}
    \item [$-$] $ \mathbf{S}^2 $ is the unit sphere of $\R^3$, whose area
    element is denoted $d\omega$;
    \item [$-$] $B$ is a positive regular function called the shock
    kernel or the cross section of collisions on which we make the
    following assumptions :

\begin{equation}\label{eq2.22}
 \;\,\left\{
               \begin{array}{ll}
               \partial^{\beta}_{\overline{p}}B\quad \mbox{is Lipschtz in the $a$ variable for
               }\; 0\leq |\beta|\leq 3\\
               \\
                 \exists C>0,\;0\leq B(a,\overline{p},\overline{q},\overline{p}',\overline{q}')\leq C &  \\
                 \\
                 (1+|\overline{p}|)^{\ell}\|\partial_{\overline{p}}^{\beta}B\|_{L^1(\R^3\times S^2)}
                 \in L^{\infty}(\R^3),\; 0\leq|\beta|\leq 3,\;0\leq  \ell \leq 3 & \\
                 \\
                 (1+|\overline{p}|)^{|\beta|-1}\partial_{\overline{p}}^{\beta}B\in L^{\infty}(\R^3\times\R^3\times S^2)
                 ,\; 1\leq|\beta|\leq 3 &
               \end{array}
               \;;
             \right.
\end{equation}
where
$\,|\overline{p}|=\left(\sum\limits_{i=1}^3(p^i)^2\right)^{1/2}$.
    \item [$-$] The conservation law $p+q=p'+q'$ splits into:
\begin{eqnarray}
\strut
       p^0+q^0&=&p'^0+q'^0 \label{eq2.23}\\
       \overline{p}+\overline{q}&=&\overline{p}'+\overline{q}'\label{eq2.24}
       \;.
      \end{eqnarray}
Equation $(\ref{eq2.23})$ is the conservation of the elementary
energy :
\begin{equation}\label{eq2.25}
    \tilde{e}=\sqrt{1+a^2(t)\big[(p^1)^2+(p^2)^2+(p^3)^2}\big]+\sqrt{1+a^2(t)\big[(q^1)^2+(q^2)^2+(q^3)^2}\big]
\end{equation}
and we interpret (\ref{eq2.24}) following R.T.Glassey \cite{RTG} by
setting :
\begin{equation}\label{eq2.26}
    \overline{p}'=\overline{p}+b(\overline{p},\overline{q},\omega)\omega
    \;;\quad
    \overline{q}'=\overline{q}-b(\overline{p},\overline{q},\omega)\omega \;; \quad \omega \in
    \mathbb{S}^2
\end{equation}
in which $b$ is a regular function of its arguments, a real valued
function which, following (\ref{eq2.8}), (\ref{eq2.23}) and
(\ref{eq2.24}) is solution of a quadratic equation whose non trivial
solution is :

\begin{equation}\label{eq2.27}
b(\overline{p},\overline{q},\omega)=\frac{2p^0q^0\tilde{e}\omega\odot(\hat{\overline{q}}-\hat{\overline{p}})}
{(\tilde{e})^2-[\omega\odot(\overline{p}+\overline{q})]^2}\;;
\end{equation}
with $\hat{\overline{p}}=\frac{\overline{p}}{p^0}$, $\tilde{e}$
given by (\ref{eq2.25}). The symbol $\odot$ is the new scalar
product defined on  $\mathbb{R}^3$ by :

\begin{equation}\label{eq2.28}
    \overline{p}\odot\overline{q}=a^2(t)[p^1q^1+p^2q^2+p^3q^3]\quad
    \mbox{and}\quad
    \|\overline{p}\|^2=a^2(t)[(p^1)^2+(p^2)^2+(p^3)^2]
    \;.
\end{equation}
Now using usual formulae, the Jacobian of the change of variables
$(\overline{p},\overline{q})\longmapsto(\overline{p}',\overline{q}')$
given by (\ref{eq2.26}) is :
\begin{equation}\label{eq2.29}
\frac{\partial(\overline{p}', \overline{q}')}{\partial(\overline{p},
\overline{q})}=-\frac{p'^0q'^0}{p^0q^0}\;.
\end{equation}
\end{enumerate}
\end{enumerate}

\begin{remark}
Our hypotheses on the shock kernel $B$ (see (\ref{eq2.22})) are the
key hypotheses in the proof of the Moser-type substitution
inequalities for the collision operator $Q$
(inequalities(\ref{eq3.2}) and (\ref{eq3.3})) which are essential in
the energy estimates method we will use. On the other hand a simple
example of function satisfying these hypotheses is given by (see
\cite{NT}):
$$
B\big(a,\overline{p}, \overline{q},
\overline{p}',\overline{q}'\big)= Ae^{-a^2-|\overline{p}|^2-
|\overline{q}|^2-|\overline{p}'|^2-|\overline{q}'|^2} \;; \quad
A>0\;.
$$
\end{remark}

\subsection{Compatibility of the equations}
\begin{enumerate}
    \item [a)] It is clear that, if $\overline{p}\longmapsto
    f(t,\overline{p})$ is invariant by $SO_3$, then the integral in
    $(\ref{eq2.9})$ is zero if $\beta=i$. In \cite{Nal}, N. Noutchegueme et al. proved that  if
$\overline{p}\longmapsto
    f(0,\overline{p})$ is invariant by $SO_3$, then for all $t>0, \, \overline{p}\longmapsto
    f(t,\overline{p})$ is invariant by $SO_3$ too. We adopt the
    assumption that $\overline{p}\longmapsto
    f(0,\overline{p})$ is invariant by $SO_3$, then by
    $(\ref{eq2.3})$
    one has
    \begin{equation}\label{eq2.30}
        \nabla_{\alpha}F^{\alpha i}=0,\; i=1,2,3\;.
    \end{equation}
Now $\R^4$ being simply connected, there exists a potential vector
$A=(A_{\mu})$ such that
$$
F_{\alpha\beta}=\nabla_{\alpha}A_{\beta}-\nabla_{\beta}A_{\alpha}=
\partial_{\alpha}A_{\beta}-\partial_{\beta}A_{\alpha}\;,
$$
 then $F_{ij}=\partial_{i}A_{j}-\partial_{j}A_{i}=0$ since
$A=A(t)$, from where we obtain $F^{ij}=F_{ij}=0$. This means that
the electromagnetic field reduces to its electric part and thus,
$(\ref{eq2.30})$ writes:
\begin{equation}\label{eq2.31}
    \partial_0F^{0i}+\Gamma_{j0}^{j}F^{0i}=0
    \;.
\end{equation}
We deduce from $(\ref{eq2.1})$ that:
\begin{equation}\label{eq2.32}
    g^{00}=g_{00}=-1;\;g_{ii}=a^2;\;g^{ii}=a^{-2};\;g_{0i}=g^{0i}=0;\;g_{ij}=g^{ij}=0\;for\;i\neq
    j
    \,.
\end{equation}
 The usual formula  $
\Gamma^\lambda_{\,\alpha\beta}=\frac{1}{2}g^{\lambda\mu}(\partial_\alpha\,
    g_{\mu\beta}+\partial_\beta\,
    g_{\alpha\mu}-\partial_\mu\,g_{\alpha\beta})$
gives, with
    $\dot{a}=\frac{da}{dt}\;:$
    \begin{equation}\label{eq2.33}
        \Gamma^0_{\;ii}=\dot{a}a;\;
      \Gamma^i_{\;i0}=\Gamma^i_{\;0i}=\frac{\dot{a}}{a};\;\Gamma^0_{\,00}=0;
      \;\Gamma^0_{\alpha\beta}=0\;\mbox{for}\;\alpha\neq\beta;\;\Gamma^k_{ij}=0
\;;
    \end{equation}
and $(\ref{eq2.31})$ is then equivalent to:
\begin{equation}\label{eq2.34}
\partial_0F^{0i}+3\frac{\dot{a}}{a}F^{0i}=0
\;.
\end{equation}

The general solution of $(\ref{eq2.34})$ is, with $a_0=a(0)$ :
\begin{equation}\label{eq2.35}
F^{0i}(t)=\left(\frac{a_0}{a}\right)^3F^{0i}(0) \;.
\end{equation}
Since $F^{ij}=F_{ij}=0$, we have directly:
\begin{equation}\label{eq2.36}
    F^{\lambda\mu}F_{\lambda\mu}=-2g_{ij}F^{0i}F^{0j};\;F_{i\lambda}F_j^{\lambda}=-2g_{ik}g_{jl}F^{0k}F^{0l};\;F_{0\lambda}F_j^{\lambda}=0;\;
    F_{0\lambda}F^{0\lambda}=g_{ij}F^{0i}F^{0j}
\end{equation}
and definition $(\ref{eq2.7})$ of $\tau_{\alpha\beta}$ gives:
\begin{equation}\label{eq2.37}
    \tau_{00}=\frac{1}{2}g_{ij}F^{0i}F^{0j};\;\tau_{0j}=0;\;\tau_{ij}=\big(\frac{1}{2}g_{ij}g_{kl}-g_{ik}g_{jl}\big)F^{0k}F^{0l}
\;.
\end{equation}

By definition $(\ref{eq2.12})$ of $H_{\alpha\beta}$, we have:
\begin{equation}\label{eq2.38}
    H_{00}=\frac{1}{2}\big(\dot{\Phi}^2+m^2\Phi^2\big);\;H_{0i}=0;\;H_{ij}=\frac{1}{2}g_{ij}\big(\dot{\Phi}^2-m^2\Phi^2\big)
\;.
\end{equation}

    \item [b)] For the Einstein equations,
    when $\alpha=0$ and $\beta=i$, we have:
    $R_{0i}-\frac{1}{2}Rg_{0i}+\Lambda g_{0i}=0$  at the l.h.s. and
    we must then have at the r.h.s:
    \begin{equation}\label{eq2.39}
T_{0i}+\tau_{0i}+K_{0i}+H_{0i}=0 \;.
\end{equation}
Since  $\overline{p}\longmapsto f(t,\overline{p})$ is invariant by
$SO_3$, we then have for $T_{\alpha\beta}$ given by $(\ref{eq2.6})$,
$T_{0i}=0$.
Equations $(\ref{eq2.37})$ and $(\ref{eq2.38})$ give
$\tau_{0i}=H_{0i}=0$ and by definition $(\ref{eq2.11})$ of
$K_{\alpha\beta}$ we will have $(\ref{eq2.39})$ if we take
$\theta_{0i}=0$.
\item [c)] For $\alpha=i,\;\beta=j$ where $i\neq j$, we have $R_{ij}-\frac{1}{2}Rg_{ij}+\Lambda
    g_{ij}=0$ at the l.h.s. We must then have at the r.h.s :
    \begin{equation}\label{eq2.40}
T_{ij}+\tau_{ij}+K_{ij}+H_{ij}=0,\;\;i\neq j \;.
\end{equation}
Equations $(\ref{eq2.38})$ gives for $i\neq j$ : $H_{ij}=0$; since
$\overline{p}\longmapsto f(t,\overline{p})$ is invariant by $SO_3$,
the definition $(\ref{eq2.6})$ of $T_{\alpha\beta}$ gives $T_{ij}=0$
if $i\neq j$. Expression $(\ref{eq2.37})$ of $\tau_{\alpha\beta}$
gives $\tau_{ij}=0$. We will then have $(\ref{eq2.40})$ if we take
$K_{ij}=-\theta_{ij}=0$, $i\neq j$.

\item [d)] Consider the case $\alpha=\beta=i\in\{1,2,3\}$:
In the Robertson-Walker space time, we have:
$g_{11}=g_{22}=g_{33}=a^2$, which implies:
$R_{11}=R_{22}=R_{33}$. So the three corresponding Einstein
equations have the same l.h.s which is, after calculations:
\begin{equation}\label{eq2.41}
\mathcal{A}=-2\ddot{a}a-(\dot{a})^2+\Lambda a^2 \;.
\end{equation}
The r.h.s of these equations must be the same. This means that, we
must have :
 \begin{equation}\label{eq2.42}
    T_{11}+\tau_{11}+K_{11}+H_{11}=T_{22}+\tau_{22}+K_{22}+H_{22}
    =T_{33}+\tau_{33}+K_{33}+H_{33}
\end{equation}
Since $\overline{p}\longmapsto f(t,\overline{p})$ is invariant by
$SO_3$ we have: $T_{11}=T_{22}=T_{33}$.
Now if in the solution $(\ref{eq2.35})$ of the equation
$(\ref{eq2.34})$ we take :
\begin{equation}\label{eq2.44}
F^{01}(0)=F^{02}(0)=F^{03}(0)
\end{equation}
 we will then have:
\begin{equation}\label{eq2.45}
F^{01}=F^{02}=F^{03} \;.
\end{equation}
We make the hypothesis $(\ref{eq2.44})$. We then have
$(\ref{eq2.45})$ and expression $(\ref{eq2.37})$ of
$\tau_{\alpha\beta}$ implies:
\begin{equation}\label{eq2.46}
\tau_{11}=\tau_{22}=\tau_{33} \;.
\end{equation}
Since $g_{ii}=a^2$, expression $(\ref{eq2.38})$ of $H_{\alpha\beta}$
implies that:
 \begin{equation}\label{eq2.47}
H_{11}=H_{22}=H_{33} \;.
\end{equation}
The relations $(\ref{eq2.42})$ then imply, given $(\ref{eq2.47})$
and since $T_{11}=T_{22}=T_{33}$, and
$K_{\alpha\beta}=-\theta_{\alpha\beta}$:
 \begin{equation}\label{eq2.48}
\theta_{11}=\theta_{22}=\theta_{33} \;.
\end{equation}
We use the hypothesis $g^{ij}\theta_{ij}=0$ (see (\ref{eq2.13})  to
have :
 \begin{equation}\label{eq2.49}
\theta_{11}+\theta_{22}+\theta_{33}=0
\end{equation}
and Equations $(\ref{eq2.48})$ and $(\ref{eq2.49})$ imply:
 \begin{equation}\label{eq2.50}
\theta_{11}=\theta_{22}=\theta_{33}=0 \;.
\end{equation}
Hence the Einstein evolution equations reduce to the single
equation:
\begin{equation}\label{eq2.51}
R_{11}-\frac{1}{2}Rg_{11}+\Lambda
g_{11}=8\pi(T_{11}+\tau_{11}+K_{11}+H_{11}) \;.
\end{equation}
\item [e)] For $\alpha = \beta=0$, we have the Hamiltonian constraint,
which writes, after calculations (see for example \cite{NT}):
 \begin{equation}\label{eq2.52}
3\Big(\frac{\dot{a}}{a}\Big)^2-\Lambda=8\pi\Big[T_{00}+\tau_{00}+K_{00}+H_{00}\Big]\;.
\end{equation}
We know that $(\ref{eq2.52})$ will be solved everywhere if and only
if it is solved for $t=0$.
 The only unknown component of
$\theta_{\alpha\beta}$ is $\theta_{00}=\theta^{00}$ which satisfies,
given $(\ref{eq2.13})$ and $u^0=1$ :
 \begin{equation}\label{eq2.53}
\dot{\theta}^{00}+3\frac{\dot{a}}{a}\theta^{00}=-\rho^2\;;
\end{equation}
whose general solution is :
\begin{equation}\label{eq2.54}
\theta^{00}=\big(\frac{a_0}{a}\big)^3\big[\theta^{00}(0)-\rho^2\big(\frac{a}{a_0}\big)^3\big].
\end{equation}
By $(\ref{eq2.54})$ ;
\begin{equation}\label{eq2.55}
\big(\theta^{00}(0)\leq 0\big)\Longrightarrow\big(\theta^{00}\leq
0\big).
\end{equation}
So we look for a solution $\theta^{00}$ of $(\ref{eq2.53})$ such
that :
 \begin{equation}\label{eq2.56}
\theta^{00}\leq 0.
\end{equation}
\end{enumerate}

\subsection{Conditions of conservation. Equation for $\Phi$}
We always have the identities :
$\nabla_{\alpha}\big(R^{\alpha\beta}-\frac{1}{2}Rg^{\alpha\beta}+\Lambda
g^{\alpha\beta}\big)=0$. We must then have at the r.h.s of equations
$(\ref{eq2.2})$ :
 \begin{equation}\label{eq2.57}
\nabla_{\alpha}T^{\alpha\beta}+\nabla_{\alpha}\tau^{\alpha\beta}+\nabla_{\alpha}K^{\alpha\beta}+\nabla_{\alpha}H^{\alpha\beta}=0
\;.
\end{equation}
 But we know from \cite{JE}, when $f$ is solution to the Boltzmann
equation, we have: $\nabla_{\alpha}T^{\alpha\beta}=0$. Then
$(\ref{eq2.57})$ reduces to:
 \begin{equation}\label{eq2.58}
\nabla_{\alpha}\tau^{\alpha\beta}+\nabla_{\alpha}K^{\alpha\beta}+\nabla_{\alpha}H^{\alpha\beta}=0
\;.
\end{equation}
 From expressions of  $\tau_{\alpha\beta}$ and $H_{\alpha\beta}$ given by
 $(\ref{eq2.7})$ and $(\ref{eq2.12})$ it is easy to see that
  \begin{numcases}
\strut
\nabla_{\alpha}\tau^{\alpha\beta}=F^{\beta}_{\,\lambda}\nabla_{\alpha}F^{\alpha\lambda} \label{eq2.59}\\
\nabla_{\alpha}H^{\alpha\beta}=\nabla^{\beta}\Phi(\Box_g\Phi-m^2\Phi)\label{eq2.60}
\end{numcases}
where $\Box_g$ is the d'Ambertian operator
$\Box_g=\nabla_{\alpha}\nabla^{\alpha}$ of $g$. But we know that
$\nabla_{\alpha}F^{\alpha 0}=0$ thus, by $(\ref{eq2.9})$ and since
$u^i=0$ we have
$$
\nabla_{\alpha}\tau^{\alpha\beta}=F^{\beta}_{\,i}\nabla_{\alpha}F^{\alpha
i}=F^{\beta}_{\,i}J^i=0 \;.
$$
We now have, given $(\ref{eq2.11})$ and $(\ref{eq2.13})$:
$\nabla_{\alpha}K^{\alpha\beta}=-\nabla_{\alpha}\theta^{\alpha\beta}=\rho^2u^{\beta}$.
So we deduce from $(\ref{eq2.58})$, using $(\ref{eq2.59})$ and
$(\ref{eq2.60})$ that:
 \begin{equation}\label{eq2.61}
 \rho^2u^{\beta}+\nabla^{\beta}\Phi(\Box_g\Phi-m^2\Phi)=0
 \;.
\end{equation}
 But for $\beta=i$, we have
 $\nabla^i\Phi=g^{i\lambda}\nabla_{\lambda}\Phi=g^{ij}\partial_i\Phi=0$
 since $\Phi=\Phi(t)$ and $u^i=0$. So from $(\ref{eq2.61})$ we have for $\beta=0$ :
  \begin{equation}\label{eq2.62}
\nabla^0\Phi(\Box_g\Phi-m^2\Phi)=-\rho^2 \;.
\end{equation}
Further,
 $\nabla^0\Phi=g^{0\lambda}\partial_{\lambda}\Phi=g^{00}\partial_0\Phi=-\dot{\Phi}$
 (since $g^{00}=-1$); and a development of $\Box_g\Phi$ gives
 $$
 \Box_g\Phi-m^2\Phi=-\big[\ddot{\Phi}+3\big(\frac{\dot{a}}{a}\big)\dot{\Phi}+m^2\Phi\big]
\;.
 $$
 From $(\ref{eq2.62})$ we then have the equation for $\Phi$ :
  \begin{equation}\label{eq2.63}
\dot{\Phi}\Big[\ddot{\Phi}+3\Big(\frac{\dot{a}}{a}\Big)\dot{\Phi}+m^2\Phi\Big]=-\rho^2
\;.
\end{equation}
Considering $(\ref{eq2.2}), (\ref{eq2.34}), (\ref{eq2.53})$ and
$(\ref{eq2.18})$ divided by $p^0$, we obtain the following system
which is equivalent to system (\ref{eq2.2})-(\ref{eq2.5}) ( we are
using $p_1=a^2p^1$):
\begin{numcases}
\strut 3\Big(\frac{\dot{a}}{a}\Big)^2-\Lambda=8\pi\Big[
a^3\int_{\mathbb{R}^3}p^0f(\overline{p})d\overline{p}+\frac{3}{2}
a^2(F^{01})^2-\theta_{00}+\frac{a^2}{2}(\dot{\Phi}^2+m^2\Phi^2)\Big]\label{eq2.64} \\
     -2\ddot{a}a-(\dot{a})+\Lambda a^2=8\pi\Big[
a^7\int_{\mathbb{R}^3}\frac{(p^1)^2}{p^0}f(\overline{p})d\overline{p}+
\frac{a^4}{2}(F^{01})^2+\frac{a^2}{2}(\dot{\Phi}^2-m^2\Phi^2)\Big]\label{eq2.65} \\
     \dot{F}^{01}+3\frac{\dot{a}}{a}F^{01}=0\label{eq2.66} \\
      \dot{\theta}^{00}+3\frac{\dot{a}}{a}\theta^{00}=-\rho^2\label{eq2.67} \\
      \dot{\Phi}\Big[\ddot{\Phi}+3\Big(\frac{\dot{a}}{a}\Big)\dot{\Phi}+m^2\Phi\Big]=-\rho^2\label{eq2.68} \\
  \frac{\partial f}{\partial
t}-2\frac{\dot{a}}{a}p^i\frac{\partial f}{\partial
p^i}-\Big(a^3F^{01}\int_{\mathbb{R}^3}f(t,\overline{q})d\overline{q}\Big)\sum\limits_{i=1}^3\frac{\partial
f}{\partial p^i}=\frac{1}{p^0}Q(f,f)\label{eq2.69}
\end{numcases}
From now we will refer to this system as the coupled
Einstein-Maxwell-Boltzmann system with cosmological constant in
presence of a massive scalar field.

\section{Local solution for the Boltzmann Equation}
In this section we recall a local existence result obtained by NN
and KNM  in \cite{NoKe} for the Boltzmann equation. The functional
spaces we will use are the same as those of that reference and they
are presented here. We will also recall the most important
substitution type inequalities concerning the collision operator $Q$
the proof of which is given in details in \cite{NoKe}. We begin with
the functional spaces.
\begin{definition}{Definition}
Let $m\in\N$, $d\in\R^+$, $T>0$,
$|\overline{p}|=\left(\sum\limits_{i=1}^3(p^i)^2\right)^{1/2}$ and
define the spaces:
\begin{enumerate}
\item
$ L^1_1(\mathbb{R}^3)=\big\{f :
\mathbb{R}^3\longrightarrow\mathbb{R}\;,\;(1+|\overline{p}|)f\in
L^1(\mathbb{R}^3)\}\; $
 and
$ \; L^2_d(\mathbb{R}^3)=\big\{f :
\mathbb{R}^3\longrightarrow\mathbb{R}\;,\;(1+|\overline{p}|)^d f\in
L^2(\mathbb{R}^3)\}.$
   \item \; $H^m_d(\mathbb{R}^3)=\big\{f :
\mathbb{R}^3\longrightarrow\mathbb{R}\,,\;
(1+|\overline{p}|)^{d+|\beta|}\partial^{\beta}_{\overline{p}}f\in
L^2(\mathbb{R}^3),\;| \beta|\leq m\big\}$.  $H^m_d(\mathbb{R}^3)$ is
a separable Hilbert space with the norm :
    $$\|f\|_{H^m_d{(\mathbb{R}^3)}}=\max\limits_{0\leq|\beta|\leq
3}\big\|(1+|\overline{p}|)^{d+|\beta|}\partial^{\beta
}_{\overline{p}}f\big\|_{L^2(\mathbb{R}^3)} \;.
$$
    \item \; $H^m_d(0, T, \mathbb{R}^3)=\big\{f :[0,
T]\times \mathbb{R}^3\longrightarrow\mathbb{R}\;,\;f\;
\mbox{\emph{continuous}},\, f(t,\cdot)\in
H^m_d(\mathbb{R}^3),\;\forall t\in[0,T]\big\}$.  Endowed with the
norm :
$$
\|f\|_{H^m_d(0, T, \mathbb{R}^3)}=\sup\limits_{t\in [0,
\,T]}\max\limits_{0\leq|\beta|\leq
3}\big\|(1+|\overline{p}|)^{d+|\beta|}\partial^{\beta
}_{\overline{p}} f(t,\cdot)\big\|_{L^2(\mathbb{R}^3)} \;;
$$
$H^m_d(0, T, \mathbb{R}^3)$ is a Banach space.
\end{enumerate}
\end{definition}
For a fixed $r>0$, we set :
\begin{equation}\label{eq3.1}
    H^m_{d,r}(0, T, \mathbb{R}^3)=\{f\in H^m_d(0, T,
\mathbb{R}^3),\;\;\|f\|_{H^m_d(0, T, \mathbb{R}^3)}\leq r\} \;.
\end{equation}
Endowed with the norm induced by $H^m_d(0, T, \mathbb{R}^3)$,\,
$H^m_{d,r}(0, T, \mathbb{R}^3)$ is a complete metric space.
\begin{remark}
We have the embeddings
\bel{inj} H^m_d( \mathbb{R}^3)\hookrightarrow
L^2_d(\mathbb{R}^3)\hookrightarrow
L^1_1(\mathbb{R}^3)\hookrightarrow L^1( \mathbb{R}^3); \;\; m\in \N,
\;d>5/2\;. \ee
\end{remark}
In \cite{NoKe} it has been proven that the functions
$(\overline{p},\overline{q},\omega)\longmapsto\partial_{\overline{p}}^{\beta}b(\overline{p},\overline{q},\omega)$,
$1\leq|\beta|\leq 3$, are bounded and therefore important results
for the collision operator $Q$ defined by $(\ref{eq2.19}),
(\ref{eq2.20})$, and $(\ref{eq2.21})$ are established. We recall
these fundamental inequalities.
\begin{proposition}
Let the collision kernel $B$ satisfies hypotheses  (\ref{eq2.22}),
and $f,g\in H^m_d(0, T, \mathbb{R}^3)$, with $d>\frac{5}{2}$. Then
$\frac{1}{p^0}Q(f,g)$ is also in $H^m_d(0, T, \mathbb{R}^3)$ and
$(f,g)\longmapsto\frac{1}{p^0}Q(f,g)$ is uniformly continuous from
$H^3_d(\R^3)\times H^3_d(\R^3)$ to $H^3_d(\R^3)$. Namely, the
following holds
    \begin{equation}\label{eq3.2}
\Big\|\frac{1}{p^0}Q(f,g)\Big\|_{H^3_{d}(\R^3)}\leq
C(T)\|f\|_{H^3_d(\R^3)}\|g\|_{H^3_d(\R^3)},\;\forall  f,g\in
H^3_{d}(\R^3) \;;
    \end{equation}
    \begin{equation}\label{eq3.3}
\Big\|\frac{1}{p^0}Q(f,f)-\frac{1}{p^0}Q(g,g)\Big\|_{H^3_{d}(\R^3)}\leq
C(T)\big(\|f\|_{H^3_d(\R^3)}+\|g\|_{H^3_d(\R^3)}\big)\|f-g\|_{H^3_d(\R^3)},\;\forall
 f,g \in H^3_{d}(\R^3) \;;
    \end{equation}
where $C(T)$ is a positive constant which only depends on $T$.
\end{proposition}
In \cite{NoKe}, NN and KNM proved, using the Faedo-Galerkin method
that the Boltzmann equation has a solution in $H^3_d(0, T,
\mathbb{R}^3)$. This is particularly important since such a solution
is of class $C^1$. We call it a regular solution. We recall here the
precise statement of this result.
\begin{proposition}  Let $f_0\in H^3_{d,r}(\R^3)$.
Then there exists $T>0$ such that the Boltzmann equation
$(\ref{eq2.69})$ has a unique solution $f$ in $H^3_d(0,T,\R^3)$ such
that $f(0,\overline{p})=f_0(\overline{p})$.
\end{proposition}
\pr See \cite{NoKe}. \qed

\section{The Einstein-Maxwell-Boltzmann system with massive scalar field as a first order system.
 The sequence iterates}

\subsection{Boltzmann Equation in covariant variables}
The Boltzmann equation simplifies when written in terms of covariant
coordinates; this is inspired by the Vlasov case where the
distribution function in these coordinates is independent of time in
the Robertson-Walker spacetime. We will then use the covariant
coordinates; of course, because of the presence of the
electromagnetic field, the Boltzmann equation do not simplifies as
in \cite{HL2, HLAR} but by using them we will get rid of the second
term in (\ref{eq2.69}). The resulting equation (\ref{eq4.30}) has
 a form to which Proposition \ref{ta} and Corollary \ref{c1}
apply (there we need our PDE to have spatially homogeneous and
bounded coefficients). Now for simplicity, instead of using lower
indices, we introduce new variables $u$ and $v$ as follows (see
\cite{HL2}). Recall that $p_i=g_{i\beta}p^{\beta}=a^2p^i$ and set
\begin{equation}\label{eq4.23}
\left\{
  \begin{array}{c}
    \overline{u}=(u^1,u^2,u^3)\,,\;\, \;\;u^i=a^2p^i\,,\;\, \;\;u^0=\sqrt{1+a^{-2}|\overline{u}|^2}=p^0 \\
   \overline{v}=(v^1,v^2,v^3)\,,\;\, \;\;v^i=a^2q^i\,,\;\, \;\;v^0=\sqrt{1+a^{-2}|\overline{v}|^2}=q^0 \\
  \end{array}
  \;.
\right.
\end{equation}
 In the change of variables $(\ref{eq2.26})$, we had the function
 $b$ defined by $(\ref{eq2.27})$, with the new scalar product
 defined by $(\ref{eq2.28})$. With the usual scalar product $(\cdot)$
 it writes :
 \begin{equation}\label{eq4.24}
b(\overline{p},\overline{q},\omega)=\frac{2a^2p^0q^0\tilde{e}\omega\cdot(\hat{\overline{q}}-\hat{\overline{p}})}
{(\tilde{e})^2-a^4[\omega\cdot(\overline{p}+\overline{q})]^2} \;.
\end{equation}
The change of variables
$$
    \overline{p}'=\overline{p}+b(\overline{p},\overline{q},\omega)\omega\,,\;\;
    \overline{q}'=\overline{q}-b(\overline{p},\overline{q},\omega)\omega \,;\; \omega \in
    \mathbb{S}^2
$$
 then writes in terms of $\overline{u}$, and
$\overline{v}$ :
$$
\left\{
  \begin{array}{c}
   p'^k=p^k+\frac{2a^2p^0q^0\tilde{e}\omega\cdot(\hat{\overline{q}}-\hat{\overline{p}})}
{(\tilde{e})^2-a^4[\omega\cdot(\overline{p}+\overline{q})]^2}w^k=\frac{u^k}{a^2}+\frac{2u^0v^0\tilde{e}\omega\cdot(\hat{\overline{v}}-\hat{\overline{u}})}
{(\tilde{e})^2-[\omega\cdot(\overline{u}+\overline{v})]^2}w^k=
\frac{1}{a^2}\Big(u^k+\frac{2a^2u^0v^0\tilde{e}\omega\cdot(\hat{\overline{v}}-\hat{\overline{u}})}
{(\tilde{e})^2-[\omega\cdot(\overline{u}+\overline{v})]^2}w^k\Big)  \\
\\
    q'^k=q^k-\frac{2a^2p^0q^0\tilde{e}\omega\cdot(\hat{\overline{q}}-\hat{\overline{p}})}
{(\tilde{e})^2-a^4[\omega\cdot(\overline{p}+\overline{q})]^2}w^k=\frac{v^k}{a^2}-\frac{2u^0v^0\tilde{e}\omega\cdot(\hat{\overline{v}}-\hat{\overline{u}})}
{(\tilde{e})^2-[\omega\cdot(\overline{u}+\overline{v})]^2}w^k=
\frac{1}{a^2}\Big(v^k-\frac{2a^2u^0v^0\tilde{e}\omega\cdot(\hat{\overline{v}}-\hat{\overline{u}})}
{(\tilde{e})^2-[\omega\cdot(\overline{u}+\overline{v})]^2}w^k\Big) \\
  \end{array}
\right.
$$
Therefore, if we set :

\begin{equation}\label{eq4.25}
u'^k=a^2p'^k,\;\;v'^k=a^2q'^k,\;\;u'^0=p'^0,\;\;v'^0=q'^0 \;;
\end{equation}
we obtain
$$
\overline{u}'=\overline{u}+\tilde{b}(\overline{u},\overline{v},\omega)\omega
\quad \mbox{and}\quad
\overline{v}'=\overline{v}-\tilde{b}(\overline{u},\overline{v},\omega)\omega\;;
$$
with
\begin{equation}\label{eq4.28}
\tilde{b}(\overline{u},\overline{v},\omega)=\frac{2a^2u^0v^0\tilde{e}\omega\cdot(\hat{\overline{v}}-\hat{\overline{u}})}
{(\tilde{e})^2-[\omega\cdot(\overline{u}+\overline{v})]^2} \;.
\end{equation}
We now write the Boltzmann equation using the variables $s, \,u$ and
$v$. Consider the change of variables :
$$
(t,\overline{p},\overline{q})\longrightarrow(s,\overline{u},\overline{v})\;\;\mbox{with}
\;
                                                               s=t\,,\;\;
                                                               \overline{u}=a^2\overline{p} \,,\;\;
                                                               \overline{v}=a^2\overline{q}
                                                             \,\,\;\;
                                                              (
                                                              \,d\overline{q}=a^{-6}d\overline{v}\;)\;.
 $$
We have
$$
    \frac{\partial f}{\partial t}=\frac{\partial s}{\partial
t}\frac{\partial \tilde{f}}{\partial s}+\frac{\partial u^i}{\partial
t}\frac{\partial \tilde{f}}{\partial u^i}=\frac{\partial
\tilde{f}}{\partial s}+2\dot{a}ap^i\frac{\partial
\tilde{f}}{\partial u^i} \quad \mbox{and} \quad
     \frac{\partial f}{\partial p^i}=\frac{\partial s}{\partial
p^i}\frac{\partial \tilde{f}}{\partial s}+\frac{\partial
u^j}{\partial p^i}\frac{\partial \tilde{f}}{\partial
u^j}=a^2\delta^j_i\frac{\partial \tilde{f}}{\partial
u^j}=a^2\frac{\partial \tilde{f}}{\partial u^i}
 $$
 where $\tilde{f}(s ,\bar{u})= f(t,\bar{p})$ and $\tilde{g}(s, \bar{v})=
 g(t,\bar{q}).$ In order to simplify our notations, we will keep $f$
 and $g$ for $\tilde{f}$ and $\tilde{g}$ respectively.
The Boltzmann equation $(\ref{eq2.69})$ in the new variables
$(s,\bar{u})$ reads
\begin{eqnarray*}
\frac{1}{u^0}\tilde{Q}(f,f)&=& \frac{\partial f}{\partial
s}+2\dot{a}ap^i\frac{\partial f}{\partial
u^i}-2\frac{\dot{a}}{a}a^2p^i\frac{\partial f}{\partial
u^i}-\left(a^3F^{01}\int_{\mathbb{R}^3}f(s,\overline{v})a^{-6}d\overline{v}\right)a^2\sum\limits_{i=1}^3\frac{\partial
f}{\partial u^i}\\
 &=&\frac{\partial f}{\partial
t}-\left(\frac{F^{01}}{a}\int_{\mathbb{R}^3}f(t,\overline{v})d\overline{v}\right)\sum\limits_{i=1}^3\frac{\partial
f}{\partial u^i} \;;
\end{eqnarray*}
with, given $(\ref{eq2.19}), (\ref{eq2.20})$ and $(\ref{eq2.21})$
$\tilde{Q}=\tilde{Q}^+-\tilde{Q}^-$ where :
\begin{equation}\label{eq4.29}
\left\{%
\begin{array}{ll}
   \tilde{Q}^+(f,
g)(t,\overline{u})=\displaystyle\int_{\mathbb{R}^{3}}\frac{a^{-3}(t)}{v^{0}}d\overline{v}\int_{S^2}f(t,\overline{u}')
g(t, \overline{v}')\tilde{B}\big(a(t),\overline{u}, \overline{v},
\overline{u}',\overline{v}',\omega\big)d\omega  & \hbox{} \\\\

    \tilde{Q}^-(f,
g)(t,\overline{u})=\displaystyle\int_{\mathbb{R}^{3}}\frac{a^{-3}(t)}{v^{0}}d\overline{v}\int_{S^2}f(t,\overline{u})
g(t, \overline{u})\tilde{B}\big(a(t),\overline{u}, \overline{v},
\overline{u}',\overline{v}',\omega\big)d\omega & \hbox{} \\
\end{array}%
\right.
\end{equation}
and where $\tilde{B}$ is defined as $B$, in terms of the new
variables $\overline{u}, \overline{v}, \overline{u}'$ and
$\overline{v}'$. Thus the Boltzmann equation  is equivalent to:

\begin{equation}\label{eq4.30}
\frac{\partial f}{\partial
t}-\left(\frac{F^{01}}{a}\int_{\mathbb{R}^3}f(t,\overline{v})d\overline{v}\right)\sum\limits_{i=1}^3\frac{\partial
f}{\partial u^i}=\frac{1}{u^0}\tilde{Q}(f,f) \;.
\end{equation}
\begin{proposition}\label{pn}
The properties established for $b$ hold for $\tilde{b}$ namely, the
functions
$(\overline{u},\overline{v},\omega)\longmapsto\partial_{\overline{u}}^{\beta}\tilde{b}
(\overline{u},\overline{v},\omega), \;1\leq|\beta|\leq 3$, are
bounded.
\end{proposition}
\pr The proof is exactly the same as for $b$  see \cite{NoKe} page
72. \qed
\begin{remark}
 An important consequence of Proposition \ref{pn} is that if we make for $\tilde{B}$ the same
assumptions  as for $B$, see (\ref{eq2.22}), we will have the
following results analogous to $(\ref{eq3.2})$ and $(\ref{eq3.3})$:
If $f,g\in H^3_d(0,T,\R^3)$, then $\frac{1}{u^0}\tilde{Q}(f,g)\in
H^3_d(0,T,\R^3)$ and the following holds
    \begin{equation}\label{eq4.31}
\Big\|\frac{1}{u^0}\tilde{Q}(f,g)\Big\|_{H^3_{d}(\R^3)}\leq
C(T)\|f\|_{H^3_d(\R^3)}\|g\|_{H^3_d(\R^3)} \;;
    \end{equation}
    \begin{equation}\label{eq4.32}
\Big\|\frac{1}{u^0}\tilde{Q}(f,f)-\frac{1}{u^0}\tilde{Q}(g,g)\Big\|_{H^3_{d}(\R^3)}\leq
C(T)\big(\|f\|_{H^3_d(\R^3)}+\|g\|_{H^3_d(\R^3)}\big)\|f-g\|_{H^3_d(\R^3)}
\;.
    \end{equation}

\end{remark}

\subsection{Change of the unknown functions in the system}
The aim of this section is write the coupled system as an equivalent
 first order system of differential equations.  The coupled system
writes, using $(\ref{eq2.64}), (\ref{eq2.65}), (\ref{eq2.66}),
(\ref{eq2.67}), (\ref{eq2.68})$ and $(\ref{eq4.30})$; the change of
variables
$(t,\overline{p},\overline{q})\longmapsto(t,\overline{u},\overline{v})$
and $d\overline{q}=a^{-6}d\overline{v}$ as: %
\begin{numcases}
{(S\,)\quad } \strut
 3\Big(\frac{\dot{a}}{a}\Big)^2-\Lambda=8\pi
a^{-3}\int_{\mathbb{R}^3}v^0f(t,\overline{v})d\overline{v}+12\pi
a^2(F^{01})^2-8\pi\theta_{00}+4\pi(\dot{\Phi}^2+m^2\Phi^2)\label{eq5.1}
\\
     -2\frac{\ddot{a}}{a}-\Big(\frac{\dot{a}}{a}\Big)^2+\Lambda=8\pi
a^{-5}\int_{\mathbb{R}^3}\frac{(v^1)^2}{v^0}f(t,\overline{v})d\overline{v}+4\pi
a^2(F^{01})^2+4\pi(\dot{\Phi}^2-m^2\Phi^2)\label{eq5.2}
\\
     \dot{F}^{01}+3\frac{\dot{a}}{a}F^{01}=0\label{eq5.3}
\\
      \dot{\theta}^{00}+3\frac{\dot{a}}{a}\theta^{00}=-\rho^2.\label{eq5.4}
      \\
      \dot{\Phi}\Big[\ddot{\Phi}+3\Big(\frac{\dot{a}}{a}\Big)\dot{\Phi}+m^2\Phi\Big]=-\rho^2\label{eq5.5}
       \\
  \frac{\partial f}{\partial
t}-\left(\frac{F^{01}}{a}\int_{\mathbb{R}^3}f(t,\overline{v})d\overline{v}\right)\sum\limits_{i=1}^3\frac{\partial
f}{\partial u^i}=\frac{1}{u^0}\tilde{Q}(f,f) \label{eq5.6}
\end{numcases}
Note that Equation $(\ref{eq5.1})$ is the Hamiltonian constraint. It
is well known that it will be satisfied everywhere if it is
satisfied for $t=0$, i.e if the initial data satisfy  the following
constraint :
\begin{equation}\label{eq5.7}
\left\{
  \begin{array}{c}
   3\Big(\frac{\dot{a}_0}{a_0}\Big)^2-\Lambda=8\pi
a_0^{-3}\int_{\mathbb{R}^3}v^0f_0(\overline{v})d\overline{v}+12\pi
a_0^2(F^{01}(0))^2-8\pi\theta_{00}(0)+4\pi(\dot{\Phi}_0^2+m^2\Phi_0^2)  \\
    with\;\; a_ 0=a(0);\; \dot{a}_0=\dot{a}(0); \;f_0(\overline{v})=f(0,\overline{v});
     \;\Phi_0=\Phi(0);\; \dot{\Phi}_0=\dot{\Phi}(0); \; F^{01}(0)= F^{01}_0  \\
  \end{array}
  \;.
\right.
\end{equation}
So we will suppose that $(\ref{eq5.7})$ holds. Therefore,
$(\ref{eq5.1})$ is solved and shall be considered as a relation
between the unknown functions. In order to have an equivalent first
 order system,
we set:
\begin{equation}\label{eq5.8}
    E=\frac{1}{a},\;U=\frac{\dot{a}}{a},\;\psi=\frac{1}{2}(\dot{\Phi})^2,
    \;\,
    Z=F^{01},\;\, W=\theta^{00}
    \;.
\end{equation}
Note that $\rho\neq 0$ implies using Equation (\ref{eq5.5}) that
$\dot{\Phi}$ does not vanish. Since $\dot{\Phi}$ is continuous, it
keeps a constant sign. We choose to look for $\dot{\Phi}$ such that
:
\begin{equation}\label{eq5.9}
    \dot{\Phi}>0
    \;.
\end{equation}
By $(\ref{eq5.9})$, $\Phi$ is increasing and we choose
$\displaystyle
    \Phi(0):=\Phi_0>0
$
 which implies
$ \displaystyle     \Phi > 0. $
We have  $\displaystyle
      \psi=\frac{1}{2}(\dot{\Phi})^2,\;\dot{\Phi}>0
      $ thus $\displaystyle\dot{\Phi}=\sqrt{2\psi};$
 $\displaystyle
 \dot{E}= -UE$
 and
        $\displaystyle\frac{\ddot{a}}{a}=
\dot{U}+U^2
. $
Therefore, we deduce from $(S)$ the equivalent system $(S')$ of
first order :
\begin{numcases}
{(S\,'\,)\quad } \strut
\dot{E}=-UE\label{eq5.12}\\
\dot{U}=-\frac{3}{2}U^2+\frac{\Lambda}{2}-4\pi
E^5\int_{\R^3}\frac{(v^1)^2}{v^0}f(t,\overline{v})d\overline{v}-
\frac{2\pi}{E^2}Z^2-2\pi(2\psi-m^2\Phi^2)
\label{eq5.13}\\
\dot{W}=-3UW-\rho^2\label{eq5.14}\\
\dot{Z}=-3UZ\label{eq5.15}\\
\dot{\Phi}=\sqrt{2\psi}\label{eq5.16}\\
\dot{\psi}=-6U\psi-m^2\Phi\sqrt{2\psi}-\rho^2\label{eq5.17}\\
\frac{\partial f}{\partial
t}-\Big(EZ\int_{\mathbb{R}^3}f(t,\overline{v})d\overline{v}\Big)\sum\limits_{i=1}^3\frac{\partial
f}{\partial u^i}=\frac{1}{u^0}\tilde{Q}(f,f)\label{eq5.18}
\end{numcases}
We will study $(S\,'\,)$ with the following initial data :
\begin{equation}\label{eq5.19}
\left\{
  \begin{array}{ll}
     E(0)=E_0=\frac{1}{a_0},\,U(0)=U_0=\frac{\dot{a}_0}{a_0},\;W(0)=W_0<0,\;Z(0)=Z_0,&  \\
     \Phi(0)=\Phi_0> 0,\;\psi(0)=\psi_0\geq
0,\;f(0,\cdot)=f_0\in H_{d,r}^3(\R^3) &\\
a_0=a(0),\;\dot{a}_0=\dot{a}(0)
  \end{array}
\right.
\end{equation}
where $H_{d,r}^3(\R^3)$ is defined by $(\ref{eq3.1})$. We suppose
that the initial data satisfy the constraint $(\ref{eq5.7})$. We
choose, given $(\ref{eq2.55})$,  $W(0)=\theta^{00}(0)<0$ and by
$(\ref{eq2.56})$, we have :
\begin{equation}\label{eq5.20}
    W\leq 0
    \;.
\end{equation}

\subsection{The sequence of iterates}
The reader may wonder why one can not solve directly Equation
 (\ref{eq5.18}) by introducing its equivalent
characteristic system. But it is not clear how to apply the
Caucy-Lipschitz Theorem to that system since the equation at hand is
an integro-differential equation in which appears the unknown $f$
and its integral. For this reason, we choose to introduce an
iterative scheme in which the characteristics method is
used to derive the solutions of the linearized equations.\\
 Let $T>0$ be given. Define on $[0,\,T]$, the functions
$E^0$, $U^0$, $W^0$, $Z^0$, $\Phi^0$, $\psi^0$ and $f^0$ by :
 $$
 E^0(t)=E_0,\; U^0(t)=U_0,\; W^0(t)=W_0,\;Z^0(t)=Z_0,\;
\Phi^0(t)=\Phi_0,\;\psi^0(t)=\psi_0\;\mbox{and}\;
f^0(t,\overline{u})=f_0(\overline{u}) \;.
$$
Now define $(E^1, U^1, W^1, Z^1, \Phi^1, \psi^1,f^1)$ as solution of
the linear system:
\begin{numcases}
{(S_0'\,)\quad}
\strut \dot{E}^1=-U^0E^0\label{eq5.21}\\
\dot{U}^1=-\frac{3}{2}(U^0)^2+\frac{\Lambda}{2}-4\pi
(E^0)^5\int_{\R^3}\frac{(v^1)^2}{v^0_0}f^0(t,\overline{v})d\overline{v}-
2\pi\frac{(Z^0)^2}{(E^0)^2}-2\pi(2\psi^0-m^2(\Phi^0)^2)
\label{eq5.22}\\
\dot{W}^1=-3U^0W^0-\rho^2\label{eq5.23}\\
\dot{Z}^1=-3U^0Z^0\label{eq5.24}\\
\dot{\Phi}^1=\sqrt{2\psi^0}\label{eq5.25}\\
\dot{\psi}^1=-6U^0\psi^0-m^2\Phi^0\sqrt{2\psi^0}-\rho^2\label{eq5.26}\\
\frac{\partial f^1}{\partial
t}-\Big(E^0Z^0\int_{\mathbb{R}^3}f^0(t,\overline{v})d\overline{v}\Big)\sum\limits_{i=1}^3\frac{\partial
f^1}{\partial u^i}=\frac{1}{u^0_0}\tilde{Q}_0\label{eq5.27}
\end{numcases}
where $u^0_0=\sqrt{1+(E^0)^2|\overline{u}|^2}$ and $\tilde{Q}_0$
stands for the collision operator $\tilde{Q}$ as defined by $E^0,
\;f^0$ and $v^0_0=\sqrt{1+(E^0)^2|\overline{v}|^2}$, with the
initial data :
$$
(E^1,\,U^1,\,W^1,\,Z^1,\,\Phi^1,\,\psi^1,\,f^1)(0)=(E_0,\,U_0,\,W_0,\,Z_0,\,\Phi_0,\,\psi_0,\,f_0)
\;.
$$
In fact by direct integrations, equations $(\ref{eq5.21})$ to
$(\ref{eq5.26})$ give $(E^1,\,U^1,\,W^1,\,Z^1,\,\Phi^1,\,\psi^1)$.
Now the partial differential equation $(\ref{eq5.27})$ is
equivalent, taking $t$ as parameter and setting
$h^1(t)=f^1(t,\overline{u}(t)),$ to its characteristic system:
 $$
 (S^0_c)\; :\; \left\{
       \begin{array}{ll}
         \frac{du^i}{dt}=-E^0Z^0\int_{\R^3}f_0(\overline{v})d\overline{v}\qquad\;\;i=1,2,3 &  \\
         \frac{dh^1}{dt}=\frac{1}{u^0_0}\tilde{Q}_0(f^0,f^0) &
       \end{array}
     \right.
     \;.
    $$
By simple integration, the characteristic system $(S^0_c)$ has a
unique solution $(\overline{u},h^1)$ which is of class $C^1$ on
$[0,\,T]$. This gives the unique solution $f^1$ of $(\ref{eq5.27})$
on $[0,\,T]$. From there one obtains existence of the solution
$(E^1,\,U^1,\,W^1,\,Z^1,\,\Phi^1,\,\psi^1,\,f^1)$ which is in
$\Big(\mathcal{C}^1\big([0,T],\R\big)\Big)^6\times\mathcal{C}^1\big([0,T]\times\R^3,
\R\big)$. %
\begin{remark}
 Given that the function $(t,\overline{u})\longmapsto
E^0Z^0\int_{\mathbb{R}^3}f^0(t,\overline{v})d\overline{v}$ is
bounded  because
$$
\Big|E^0Z^0\int_{\mathbb{R}^3}f^0(t,\overline{v})d\overline{v}\Big|\leq
E_0Z_0\|f_0\|_{H^3_{d,r}}(\R^3)\leq E_0Z_0r \;;
$$
and that its derivatives with respect to $\overline{u}$ are zero, by
Corollary \ref{c1}, Inequality $(\ref{eq4.22})$ implies that $f^1$
satisfies:
 \begin{equation}\label{eq5.28}
   e^{-\delta_1t} \|f^1(t,\cdot)\|^2_{H^3_d(\R^3)}\leq
 \|f_0\|^2_{H^3_d(\R^3)}+C_1\int_0^te^{-\delta_1s}\Big\|\frac{1}{u^0_0}\tilde{Q}_0(s,\cdot)\Big\|^2_{H^3_d(\R^3)}ds\;;
 \end{equation}
 where $\delta_1$ and  $C_1$ are positive constants which depends $E_0, |Z_0|$ and
 $r$.
 We deduce from $(\ref{eq5.28})$ given $(\ref{eq4.31})$ that:
 $$
 e^{-\delta_1t}\|f^1(t,\cdot)\|^2_{H^3_d(\R^3)}\leq
 \|f_0\|^2_{H^3_d(\R^3)}+C_2\int_0^t e^{-\delta_1s}\big\|f^0(s,\cdot)\big\|^4_{H^3_d(\R^3)}ds
 \;.
 $$
 Here $C_2$ is a positive constant which depends on $E_0,\, |Z_0|, \, r$ and
 $T$ and  since $f^0=f_0\in H^3_{d,r}(\R^3)$ we deduce  that the solution $f^1$ of
$(\ref{eq5.27})$ is in $C\big([0,T],H^3_d(\R^3)\big)$.
\end{remark}
Now we can iterate as follows. Suppose that the set of functions
$(E^n,U^n,W^n,Z^n,\Phi^n,\psi^n,f^n)$ are given in
$\Big(\mathcal{C}^1([0,T],\R)\Big)^6\times\mathcal{C}^1\big([0,T],H^3_d(\R^3)\big)$
and define
$(E^{n+1},U^{n+1},W^{n+1},Z^{n+1},\Phi^{n+1},\psi^{n+1},f^{n+1})$ as
the solution of the linear system:
 \begin{numcases}
 {(S_n'\,)\;}
\strut
\dot{E}^{n+1}=-E^nU^n\label{eq5.30}\\
\dot{U}^{n+1}=-\frac{3}{2}(U^n)^2+\frac{\Lambda}{2}-4\pi
(E^n)^5\int_{\R^3}\frac{(v^1)^2}{v^0_n}f^n(t,\overline{v})d\overline{v}
-2\pi\frac{(Z^n)^2}{(E^n)^2}-2\pi(2\psi^n-m^2(\Phi^n)^2)\label{eq5.31}\\
\dot{W}^{n+1}=-3U^nW^n-\rho^2\label{eq5.32}\\
\dot{Z}^{n+1}=-3U^nZ^n\label{eq5.33}\\
\dot{\Phi}^{n+1}=\sqrt{2\psi^n}\label{eq5.34}\\
\dot{\psi}^{n+1}=-6U^n\psi^n-m^2\Phi^n\sqrt{2\psi^n}-\rho^2\label{eq5.35}\\
\frac{\partial f^{n+1}}{\partial
t}-\Big(E^nZ^n\int_{\mathbb{R}^3}f^n(t,\overline{v})d\overline{v}\Big)\sum\limits_{i=1}^3\frac{\partial
f^{n+1}}{\partial u^i}=\frac{1}{u^0_n}\tilde{Q}_n\label{eq5.36}
\end{numcases}
 where $u^0_n=\sqrt{1+|E^n|^2|\overline{u}|^2}$ and $\tilde{Q}_n$
stands for the collision operator $\tilde{Q}$ as defined by $E^n,
\;f^n$ and $v^0_n=\sqrt{1+(E^n)^2|\overline{v}|^2}$,
 with initial data:
 $$
 (E^{n+1},\,U^{n+1},\,W^{n+1},\,Z^{n+1},\,\Phi^{n+1},\,\psi^{n+1},\,f^{n+1})(0)=(E_0,\,U_0,\,W_0,\,Z_0,
 \,\Phi_0,\,\psi_0,\,f_0)
 \;.
 $$
We use the method of characteristics as we did before to obtain that
system $(S_n')$ has a unique solution
$(E^{n+1},U^{n+1},W^{n+1},Z^{n+1},\Phi^{n+1},\psi^{n+1},f^{n+1})$ in
$\Big(\mathcal{C}^1([0,T],\R)\Big)^6\times\mathcal{C}^1\big([0,T];H^3_d(\R^3)
\big)$.
 We have thus constructed a sequence
 $(E^n,U^n,W^n,Z^n,\Phi^n,\psi^n,f^n)\in\big(\mathcal{C}^1([0,T]),\R\big)^6\times\mathcal{C}^1\big([0,T],H^3_d(\R^3)\big)$
 defined on $[0, T]$ for an arbitrary $T>0$  and we want to show
 that this sequence converges
 to the solution $(E,U,W,Z,\Phi,\psi,f)$ of the system $(S')$ on a
 maximal interval $[0, T_*)$ where $T_*$ is to be defined.

\section{Local existence of solutions of the Einstein-Maxwell-Boltzmann
system with massive scalar field} \label{s6} In this Section, by
using the techniques of energy estimates (see section \ref{Apend})
we derive a local existence theorem (Theorem \ref{MT}) to the
Einstein-Maxwell-Boltzmann system with massive scalar field as the
limit in our functional spaces  of the sequence we have just
constructed.

\subsection{ Boundedness properties of the iterates}
We start this section with some boundedness properties of the
sequence of the previous section. We have the following
\begin{proposition}\label{p1.1}
Let $f_0\in H_{d,r}^3(\R^3)$. There exists $T_0>0$, independent of
$n$, such that, the sequence $\left(X^n\right)$ where
 $X^n=(E^n,U^n,W^n,Z^n,\Phi^n,\psi^n,f^n)$  is uniformly bounded on
$[0,\,T_0]$.
\end{proposition}

\pr Set
$$\|X^n(t)\|=|E^n(t)|+|U^n(t)|+|W^n(t)|+|Z^n(t)|+|\Phi^n(t)|+|\psi^n(t)|+\|f^n(t,\cdot)\|_{H^3_d(\R^3)}$$
and
$$C_0= |E^0|+|U^0|+|W^0|+|Z^0|+|\Phi^0|+|\psi^0|+r
\;.
$$
We will prove by induction that, there exists $T_0>0$ such that
$\|X^n(t)\|\leq 2C_0, \; \;\forall n\in\N, \;\;\forall t\in [0,\,T_0]$. \\
 $\bullet$ For $n=0$, we have $\|X^0\|\leq C_0\leq 2C_0$.\\
 $\bullet$ Let $n\in\N$ and suppose that $\forall k\,\leq n$, we have
 $\|X^k(t)\|\leq 2C_0$.
 We want to show that $\|X^{n+1}(t)\|\leq
 2C_0$ for all $t\in [0,T_0]$, the choice of $T_0$ will be given shortly.
 \begin{enumerate}
    \item [-] Integrate equation $(\ref{eq5.30})$ on $[0,\,t]$, we
    have :
    $$
E^{n+1}(t)=E_0-\int_0^tE^n(s)U^n(s)ds \;;
$$
then
\begin{center}
\qquad\qquad\qquad\qquad\qquad\qquad\qquad\qquad $|E^{n+1}(t)|\leq
|E_0|+A_1t$\;;\hfill{(a)}
\end{center}
where $A_1>0$ is a constant which only depends on $C_0$.
    \item [-] Integrate equation $(\ref{eq5.31})$ on $[0,\,t]$, we
    have
\begin{equation}\label{eq6.1}
U^{n+1}(t)=U_0+
\int_0^t\Big[-\frac{3}{2}(U^n)^2+\frac{\Lambda}{2}-4\pi
(E^n)^5\int_{\R^3}\frac{(v^1)^2}{v^0_n}f^n(s,\overline{v})d\overline{v}
-2\pi\frac{(Z^n)^2}{(E^n)^2}-2\pi(2\psi^n-m^2(\phi^n)^2)\Big]ds .
\end{equation}
We had, by $(\ref{eq5.30})$ : $\dot{E}^n=-U^{n-1}E^{n-1}$ which
implies : $|\dot{E}^n|\leq 4C_0^2$ i.e. $-4C_0^2\leq \dot{E}^n\leq
4C_0^2$.\\
Integrating this last inequality gives $E_0-4C_0^2t\leq E^n(t)$. But
$E_0=\frac{1}{a_0}>0$. So if we take $t$ such that
$0<4C_0^2t<\frac{E_0}{2}$, then $E_0-4C_0^2t\geq \frac{E_0}{2}$;
then $\frac{E_0}{2}\leq E^n(t)$ i.e.
$\frac{1}{E^n(t)}\leq\frac{2}{E_0}$.
This proves that one can find two constants $t_1>0$ and $A_2>0$ such
that, identity $(\ref{eq6.1})$ gives the following:
\begin{center}
\qquad\qquad\qquad\qquad\qquad\qquad\qquad\qquad $ |U^{n+1}(t)|\leq
|U_0|+A_2t\;, \forall\, 0\leq t\leq t_1\;. $\hfill{(b)}
\end{center}

    \item [-] Integrate equation $(\ref{eq5.32})$ on $[0,\,t]$, we
    have: $
W^{n+1}(t)=W_0-\int_0^t(3U^n(s)W^n(s)+\rho^2)ds $. Then there exists
a constant $A_3>0$ such that :
\begin{center}
\qquad\qquad\qquad\qquad\qquad\qquad\qquad\qquad $ |W^{n+1}(t)|\leq
|W_0|+A_3t\;. $\hfill{(c)}
\end{center}
    \item [-] Integrate equation $(\ref{eq5.33})$ on $[0,\,t]$, we
    have: $
Z^{n+1}(t)=Z_0-3\int_0^tU^n(s)Z^n(s)ds $. Then there exists a
constant $A_4>0$ such that :
\begin{center}
\qquad\qquad\qquad\qquad\qquad\qquad\qquad\qquad $ |Z^{n+1}(t)|\leq
|Z_0|+A_4t \;.$\hfill{(d)}
\end{center}
    \item [-] Integrate equation $(\ref{eq5.34})$ on $[0,\,t]$, we
    have : $
\Phi^{n+1}(t)=\Phi_0+\int_0^t\sqrt{2\psi^n(s)}ds .$ This shows that
there exists a constant $A_5>0$ such that :
\begin{center}
\qquad\qquad\qquad\qquad\qquad\qquad\qquad\qquad $
|\Phi^{n+1}(t)|\leq |\Phi_0|+A_5t  \;.$\hfill{(e)}
\end{center}
\item [-] Integrate equation $(\ref{eq5.35})$ on $[0,\,t]$, we
    have : $$
\psi^{n+1}(t)=\psi_0-\int_0^t\big(6U^n(s)\psi^n(s)+m^2\Phi^n(s)\sqrt{2\psi^n(s)}+\rho^2\big)ds
\;.$$
 This shows that there exists a constant $A_6>0$ such that :
\begin{center}
\qquad\qquad\qquad\qquad\qquad\qquad\qquad\qquad $
|\psi^{n+1}(t)|\leq |\psi_0|+A_6t  \;.$\hfill{(f)}
\end{center}
    \item [-] Now we want to use Corollary \ref{c1} to obtain a
    bound for $\|f^{n+1}(t,\cdot)\|^2_{H^3_d(\R^3)}$. Observe that
    Equation (\ref{eq5.36}) is of the form (\ref{eq4.1}) with $a= a(t)=
    E^nZ^n\int_{\mathbb{R}^3}f^n(t,\overline{v})d\overline{v}$,
    $\,b\equiv 0$
    and since by induction hypothesis
$$
\Big| E^nZ^n\int_{\mathbb{R}^3}f^n(t,\overline{v})d\overline{v}
\Big|\leq 4C_0^2\|f^n(t, \cdot)\|_{L^1(\R^3)}\leq  C(C_0)\|f^n(t,
\cdot)\|_{H^3_d(\R^3)}\leq  C(C_0) \;;
$$
one can use inequality (\ref{eq4.22}) of Corollary \ref{c1} (with
$\frac{1}{\kappa}=C(C_0), \; k=n=3$) to obtain

    $$
\|f^{n+1}(t,\cdot)\|^2_{H^3_d(\R^3)}\leq
\Big(\|f_0\|^2_{H^3_d(\R^3)}+C(C_0)\int_0^t\|f^n(s,\cdot)\big\|^2_{H^3_d(\R^3)}ds\Big)e^{C(C_0)t}
\;.
$$
This proves that :
$$
\|f^{n+1}(t,\cdot)\|^2_{H^3_d(\R^3)}\leq
e^{C(C_0)t}\left(\|f_0\|^2_{H^3_d(\R^3)}+ C(C_0)t\right) \;.
$$
We take $t_2>0$ such that $e^{C(C_0)t}\leq 1, \,\forall
t\in[0,\,t_2]$. Then, there exists a constant $A_7>0$, such that for
$ t\in[0,\,t_2]$:
\begin{center}
\qquad\qquad\qquad\qquad\qquad\qquad $
\|f^{n+1}(t)\|_{H^3_d(\R^3)}\leq
 \|f_0\|_{H^3_d(\R^3)}+A_7\sqrt{t} \;.$   \hfill{(g)}
\end{center}
 \end{enumerate}
 Now we add inequalities $(a), (b), (c), (d), (e), (f), (g)$ to
 obtain, for $t\leq \min(t_1,t_2)$ :
 $$
\|X^{n+1}(t,\cdot)\|\leq
 C_0+\big(\sum\limits_{i=1}^7A_i\big)(t+\sqrt{t})
 \;.
$$
Now choose $t_3>0$ such that for $\forall\, 0\leq t\leq t_3$ we have
$\big(\sum\limits_{i=1}^7A_i\big)(t+\sqrt{t})\leq C_0$. Finally, by
setting $T_0=\min(t_1,t_2,t_3)$ we obtain that
$$
\|X^{n+1}(t,\cdot)\|\leq 2C_0\qquad \mbox{for  all}\quad 0\leq t\leq
T_0\;;
$$
and the proof is complete.
\qed\\
We have the following
\begin{lemma} \label{lll1}
  Let $n$ be a non negative integer, then for all $t\in [0,T_0]$; we have the following
  inequalities
  \be\label{N1}
  \left||E^n|^5\int_{\R^3}\frac{|v^1|^2}{v^0_n}f^n(\overline{v})d\overline{v}-
  |E^{n-1}|^5\int_{\R^3}\frac{|v^1|^2}{v^0_{n-1}}f^{n-1}(\overline{v})d\overline{v}\right|
  \leq C(C_0)\left(|E^n-E^{n-1}|+\|f^n-f^{n-1}\|_{H^2_d(\R^3)}\right)
\ee and
  \be\label{N2}
  \left\|\frac{1}{u^0_n}\tilde{Q}_n-\frac{1}{u^0_{n-1}}\tilde{Q}_{n-1}\right\|_{H^2_d(\R^3)}\leq
C(C_0)\left(|E^n-E^{n-1}|+\|f^n-f^{n-1}\|_{H^2_d(\R^3)}\right)
\;.
  \ee
\end{lemma}
\pr First, we prove (\ref{N1}).
 We have:
 \bea
\lefteqn{\left||E^n|^5\int_{\R^3}\frac{|v^1|^2}{v^0_n}f^n(\overline{v})d\overline{v}-
  |E^{n-1}|^5\int_{\R^3}\frac{|v^1|^2}{v^0_{n-1}}f^{n-1}(\overline{v})d\overline{v} \right|}
  \\
  &=&\left|\left(|E^n|^5-|E^{n-1}|^5\right)\int_{\R^3}\frac{|v^1|^2}{v^0_n}f^n(\overline{v})d\overline{v}+
  |E^{n-1}|^5\left(\int_{\R^3}\frac{|v^1|^2}{v^0_n}f^n(\overline{v})-\frac{|v^1|^2}
  {v^0_{n-1}}f^{n-1}(\overline{v})\right)d\overline{v}\right|\\
  &\leq& C(C_0)|E^n-E^{n-1}|+|E^{n-1}|^5\left|\int_{\R^3}\Big(\frac{1}{v^0_{n}}-\frac{1}{v^0_{n-1}}
  \Big)|v^1|^2f^n(\overline{v})d\overline{v}
  +\int_{\R^3}\frac{|v^1|^2}{v^0_{n-1}}\big(f^n-f^{n-1}\big)(\overline{v})d\overline{v}\right|\\
  &\leq& C(C_0)\left(|E^n-E^{n-1}|+\|f^n-f^{n-1}\|_{H^2_d(\R^3)}\right)
\;.
  \eea
Secondly we prove (\ref{N2}).
From now, we choose to write $\tilde{Q}(E^n,f^{n},f^{n})$ for
$\tilde{Q}_n$ to take advantage of the presence of $E^n$ in
$\tilde{Q}_n$.
  We have
\bean\label{N3}
\frac{1}{u^0_n}\tilde{Q}_n-\frac{1}{u^0_{n-1}}\tilde{Q}_{n-1}&=&\frac{1}{u^0_{n}}
\left(\tilde{Q}(E^n,f^n,f^n)-\tilde{Q}(E^{n-1},f^n,f^n)\right)+
\left(\frac{1}{u^0_n}-\frac{1}{u^0_{n-1}}\right)\tilde{Q}(E^{n-1},f^n,f^n)\nonumber\\
&&
 +\frac{1}{u^0_{n-1}}\left(\tilde{Q}(E^{n-1},f^n,f^n)-\tilde{Q}(E^{n-1},f^{n-1},f^{n-1})\right)
\;.
  \eean
As far as the first term of (\ref{N3}) is concerned, we have:
 \bea
\frac{1}{u^0_{n}}\left(\tilde{Q}(E^n,f^{n},f^{n})-\tilde{Q}(E^{n-1},f^{n},f^{n})\right)
&=&
\frac{1}{u^0_{n}}\left[\tilde{Q}^+(E^n,f^{n},f^{n})-\tilde{Q}^+(E^{n-1},f^{n},f^{n})\right]
\\
&&-
\frac{1}{u^0_{n}}\left[\tilde{Q}^-(E^n,f^{n},f^{n})-\tilde{Q}^-(E^{n-1},f^{n},f^{n})\right]
\\
&=:& (I)+(II) \;. \eea
Let us estimate the first term. We have
\bea (I)
   &=& \frac{1}{u^0_{n}}\bigg[\int_{\mathbb{R}^{3}}\frac{(E^n)^3}{v_n^{0}}d\overline{v}\int_{S^2}f^n(t,\overline{u}')
f^n(t,\overline{v}')\tilde{B}(E^n,\overline{u},\overline{v},\overline{u}',\overline{v}')d\omega \\
 &&\qquad\qquad\qquad\qquad\qquad\qquad -
\int_{\R^3}\frac{(E^{n-1})^3}{v_{n-1}^0}d\overline{v}\int_{S^2}f^n(t,\overline{u}')f^n(t,\overline{v}')
\tilde{B}(E^{n-1},\overline{u},\overline{v},\overline{u}',\overline{v}')
d\omega\bigg] \\
   &=&
\frac{1}{u^0_{n}}\bigg[\int_{\mathbb{R}^{3}}\frac{(E^n)^3}{v_n^{0}}d\overline{v}\int_{S^2}f^n(t,\overline{u}')
f^n(t,\overline{v}')\big(\tilde{B}(E^n,\overline{u},\overline{v},\overline{u}',\overline{v}')-
\tilde{B}(E^{n-1},\overline{u},\overline{v},\overline{u}',\overline{v}')\big)d\omega \\
 &&+
\int_{\R^3}\left(\frac{(E^n)^3}{v_n^0}-\frac{(E^{n-1})^3}{v_{n-1}^0}\right)d\overline{v}\int_{S^2}f^n(t,\overline{u}')f^n(t,\overline{v}')
\tilde{B}(E^{n-1},\overline{u},\overline{v},\overline{u}',\overline{v}')
d\omega\bigg]
\\
&=& (I_1)+ (I_2)
\;. \eea
Note that
\begin{eqnarray*}
 \frac{(E^n)^3}{v_n^0}-\frac{(E^{n-1})^3}{v_{n-1}^0} &=& \frac{(E^n)^3\big(v^0_{n-1}-v^0_n\big)+v^0_n\big((E^{n})^3-(E^{n-1})^3\big)}{v^0_nv^0_{n-1}} \\
   &=& \frac{(E^{n})^3-(E^{n-1})^3}{v^0_{n-1}}+(E^n)^3\frac{\sqrt{1+(E^{n-1})^2|\overline{v}|^2}-\sqrt{1+(E^n)^2|\overline{v}|^2}}{v^0_nv^0_{n-1}} \\
   &=&\frac{E^n-E^{n-1}}{v^0_{n-1}}\underbrace{\left[(E^n)^2+E^nE^{n-1}+(E^{n-1})^2+\frac{(E^n)^3(E^n+E^{n-1})
   |\overline{v}|^2}{v^0_n(v^0_n+v^0_{n-1})}\right]}_{\; :\,=\, \xi^n}
   \;;
\end{eqnarray*}
thus
\bel{er} I_2= \left(E^n-E^{n-1}\right)\frac{u^0_{n-1}}{u^0_n}\,
\frac{1}{u^0_{n-1}}\int_{\R^3}\frac{\xi^n
\,d\bar{v}}{v^0_{n-1}}\int_{S^2}f^n(t,\overline{u}')f^n(t,\overline{v}')
\tilde{B}(E^{n-1},\overline{u},\overline{v},\overline{u}',\overline{v}')
d\omega \;. \ee
Before continuing, we point out the following Lemma.
\begin{lemma}
\label{la1}
  If $f$ and $g$ are functions such that the partial derivatives of
  $ f$ up to order $k$ are uniformly bounded then
  $$
\|fg\|_{H^k_{d}(\R^3)}\leq C \|g\|_{H^k_{d}(\R^3)}
  $$
  where the positive constant $C$ only depends on the bounds of $f$
  and its derivatives.
\end{lemma}
This Lemma shows that in order to control $H^2_{d}-$norm of $(I_2)$,
we then need to show that the function $\bar{u}\mapsto
\frac{u^0_{n-1}}{u^0_n}$ and its derivatives up to order two are
uniformly bounded. Note that for all $n\in \N, \;0<
\frac{E_0}{2}\leq E^n\leq 2C_0$ thus we have
$$
\frac{u^0_{n-1}}{u^0_n}=
\frac{\sqrt{1+(E^{n-1})^2|\bar{u}|^2}}{\sqrt{1+(E^{n})^2|\bar{u}|^2}}\leq
C(C_0,E_0) \;;
$$
$$
\left|\partial_i\left(
\frac{u^0_{n-1}}{u^0_n}\right)\right|=\left|\frac{(E^{n-1})^2
u^i}{u^0_{n-1}u^0_{n}}-\frac{(E^{n})^2
u^iu^0_{n-1}}{(u^0_{n})^3}\right|\leq C(C_0,E_0) \;;
$$
and
\begin{eqnarray}\label{q2}
\left|\partial^2_{ij}\left(\frac{u^0_{n-1}}{u^0_n}\right)\right|&=&\left|\frac{(E^{n-1})^2\delta_j^i}{u^0_nu^0_{n-1}}
-2\frac{(E^n)^2(E^{n-1})^2u^iu^j}{(u^0_n)^3u^0_{n-1}}-\frac{(E^{n-1})^4u^iu^j}{u^0_n(u^0_{n-1})^3}
-\frac{\delta_j^i(E^n)^2u^0_{n-1}}{(u^0_n)^3}
+3\frac{(E^n)^4u^iu^ju^0_{n-1}}{(u^0_n)^5}\right| \nonumber
\\
&\leq& C(C_0,E_0) \;.
\end{eqnarray}
Lemma \ref{la1} implies that
$$
\|(I_2)\|_{H^2_d(\R^3)}\leq C(C_0)|E^n-E^{n-1}|\left\|
\frac{1}{u^0_{n-1}}\int_{\R^3}\frac{\xi^n
\,d\bar{v}}{v^0_{n-1}}\int_{S^2}f^n(t,\overline{u}')f^n(t,\overline{v}')
\tilde{B}(E^{n-1},\overline{u},\overline{v},\overline{u}',\overline{v}')
d\omega  \right\|_{H^2_d(\R^3)} \;.
$$
Since $\xi^n$ is bounded and does not depend on $\bar{u}$, the
$H^2_d(\R^3)-$norm at the r.h.s. of the previous inequality will
give an estimate similar to (\ref{eq4.31}). More precisely, we have
\bel{a11} \|(I_2)\|_{H^2_d(\R^3)}\leq
C|E^n-E^{n-1}|\,\|f^n\|^2_{H^2_d(\R^3)}\leq  C|E^n-E^{n-1}| \;. \ee
In order to obtain an estimate for the $H^2_d(\R^3)$-norm of the
term $(I_1)$, we proceed exactly as in the proof of  Proposition 3.6
of \cite{NoKe}, page 88. The term $B$  in that reference is replaced
by the difference
$\big(\tilde{B}(E^n,\overline{u},\overline{v},\overline{u}',\overline{v}')-
\tilde{B}(E^{n-1},\overline{u},\overline{v},\overline{u}',\overline{v}')\big)$
and we use instead the property that $B$ and its derivatives are
Lipschitz continuous. This leads to an estimate of the form
\bel{a12} \|(I_1)\|_{H^2_d(\R^3)}\leq
C|E^n-E^{n-1}|\,\|f^n\|^2_{H^2_d(\R^3)}\leq  C|E^n-E^{n-1}| \;; \ee
and then
 \bea
\left\|\frac{1}{u^0_{n}}\left[\tilde{Q}^+(E^n,f^{n},f^{n})-\tilde{Q}^+(E^{n-1},f^{n},f^{n})\right]\right\|_{H^2_d(\R^3)}\leq
C|E^n-E^{n-1}| \;.
 \eea
Similarly we have:
 \bea
\left\|(II)\right\|_{H^2_d(\R^3)}=
\left\|\frac{1}{u^0_{n}}\left[\tilde{Q}^-(E^n,f^{n},f^{n})-\tilde{Q}^-(E^{n-1},f^{n},f^{n})\right]\right\|_{H^2_d(\R^3)}\leq
C|E^n-E^{n-1}| \;;
 \eea
 and we deduce that:
 \be\label{N4}
\left\|\frac{1}{u^0_{n}}\left[\tilde{Q}(E^n,f^{n},f^{n})-\tilde{Q}(E^{n-1},f^{n},f^{n})\right]\right\|_{H^2_d(\R^3)}\leq
C|E^n-E^{n-1}| \;.
 \ee
As far as the second term of (\ref{N3}) is concerned, we have:
 \bean
\left\|\left(\frac{1}{u^0_n}-\frac{1}{u^0_{n-1}}\right)\tilde{Q}(E^{n-1},f^n,f^n)\right\|_{H^2_d(\R^3)}&=&
\left\|\frac{u^0_{n-1}-u^0_n}{u^0_n}\,\frac{1}{u^0_{n-1}}\tilde{Q}(E^{n-1},f^n,f^n)\right\|_{H^2_d(\R^3)}
\nonumber
\\
&=&
\big|(E^{n-1})^2-(E^{n})^2\big|\left\|\frac{|\bar{u}|^2}{(u^0_{n-1}+u^0_n)u^0_n}\,\frac{1}{u^0_{n-1}}\tilde{Q}(E^{n-1},f^n,f^n)
\right\|_{H^2_d(\R^3)} \nonumber
\\
&\leq&
C\big|E^{n-1}-E^{n}\big|\left\|\frac{1}{u^0_{n-1}}\tilde{Q}(E^{n-1},f^n,f^n)
\right\|_{H^2_d(\R^3)}
\nonumber
\\
&\leq& C\big|E^{n-1}-E^{n}\big|\left\|f^n \right\|^2_{H^2_d(\R^3)}
\nonumber \\
&\leq& C\big|E^{n-1}-E^{n}\big| \;; \eean
where in the first inequality we have used the fact that the
function $\bar{u}\mapsto \frac{|\bar{u}|^2}{(u^0_{n-1}+u^0_n)u^0_n}$
and its derivatives up to order two are bounded (the details of
computations can be found in\cite{KNM}).
 For the third term of (\ref{N3}), applying directly (\ref{eq4.32})
gives
 \be\label{N6}
\left\|\frac{1}{u^0_{n-1}}\left[\tilde{Q}(E^{n-1},f^n,f^n)-\tilde{Q}(E^{n-1},f^{n-1},f^{n-1})\right]\right\|_{H^2_d(\R^3)}\leq
C\|f^n-f^{n-1}\|_{H^2_d(\R^3)} \;.
 \ee
Finally adding (\ref{N4})-(\ref{N6}) gives (\ref{N2}) and the proof
is complete.
\qed
\begin{proposition} \label{p6.1} The hypotheses are those of the previous
Proposition. Set  $Y^n=(E^n,U^n,W^n,Z^n,\Phi^n,\psi^n)$ then the
sequences $(X^n), \; (Y^n)$ and $(\frac{\partial f^n}{\partial t})$
are Cauchy sequences respectively in the Banach spaces
$\big(C^0([0,T_0];\R)\big)^6\times
C^0\big([0,T_0];H^2_d(\R^3)\big)$,
$(C^1([0,T_0];\R))^6$ and $C^0\big([0,T_0];H^1_d(\R^3)\big)$
possibly for smaller $T_0$.
\end{proposition}
\pr In what follows, the constant $C$ only depends on $C_0$ and
$T_0$ and may be different from line to line.
\begin{enumerate}
\item We first prove that $(X^n)$ is a Cauchy sequence in
$\big(C^0([0,T_0];\R)\big)^6\times
C^0\big([0,T_0];H^2_d(\R^3)\big)$.
\begin{enumerate}
   \item [$\bullet$] We integrate Equation (\ref{eq5.30}) and obtain:
\begin{equation}\label{eq6.2}
|E^{n+1}(t)-E^n(t)|^2\leq
C\int_0^t\big(|E^n(s)-E^{n-1}(s)|^2+|U^n(s)-U^{n-1}(s)|^2\big)ds \;.
\end{equation}
    \item [$\bullet$] From Equation (\ref{eq5.31}) we have:
\begin{eqnarray*}
  U^{n+1}(t)-U^n(t) &=& -\displaystyle\int_0^t\bigg[\frac{3}{2}\big((U^n(s))^2-(U^{n-1}(s))^2\big)+
2\pi\Big(\frac{(Z^n)^2(s)}{(E^n)^2(s)}-\frac{(Z^{n-1})^2(s)
}{(E^{n-1})^2(s)}\Big)\\
   && +\; 4\pi
\Big((E^{n})^5(s)\int_{\R^3}\frac{|v^1|^2}{v^0_n}f^n(s,\overline{v})d\overline{v}-
(E^{n-1})^5(s)\int_{\R^3}\frac{|v^1|^2}{v^0_{n-1}}f^{n-1}(s,\overline{v})d\overline{v}\Big) \\
   &&+ \;4\pi\big(\psi^n(s)-\psi^{n-1}(s)\big)-2\pi
m^2\big((\Phi^n(s))^2-(\Phi^{n-1}(s))^2\big)\bigg]ds\;;
\end{eqnarray*}
which implies (since (\ref{N1})) that
\begin{eqnarray*}
\lefteqn{
 |U^{n+1}(t)-U^n(t)|\leq
C \int_0^t\Big[|E^n(s)-E^{n-1}(s)|+|U^n(s)-U^{n-1}(s)|+|Z^n(s)-Z^{n-1}(s)|}\\
 &&
+|\Phi^n(s)-\Phi^{n-1}(s)|+|\psi^n(s)-\psi^{n-1}(s)|
+\|f^n(s,\cdot)-f^{n-1}(s,\cdot)\|_{H^2_d(\R^3)}\Big]ds \nonumber
\;.
\end{eqnarray*}
Thus,
\begin{eqnarray}\label{eq6.3}
\lefteqn{
 |U^{n+1}(t)-U^n(t)|^2\leq
C\int_0^t\Big[|E^n(s)-E^{n-1}(s)|^2+|U^n(s)-U^{n-1}(s)|^2+|Z^n(s)-Z^{n-1}(s)|^2}\nonumber\\
 &&
+|\Phi^n(s)-\Phi^{n-1}(s)|^2+|\psi^n(s)-\psi^{n-1}(s)|^2
+\|f^n(s,\cdot)-f^{n-1}(s,\cdot)\|^2_{H^2_d(\R^3)}\Big]ds\;.
\end{eqnarray}
    \item [$\bullet$] Integrating Equation (\ref{eq5.32}) gives:
\begin{equation}\label{eq6.4}
|W^{n+1}(t)-W^n(t)|^2\leq
C\int_0^t\big(|U^n(s)-U^{n-1}(s)|^2+|W^n(s)-W^{n-1}(s)|^2\big)ds \;.
\end{equation}

    \item[$\bullet$] Similarly, Equation (\ref{eq5.33}) gives:
\begin{equation}\label{eq6.5}
|Z^{n+1}(t)-Z^n(t)|^2\leq C
\int_0^t\big(|U^n(s)-U^{n-1}(s)|^2+|Z^n(s)-Z^{n-1}(s)|^2\big)ds \;.
\end{equation}
 %
 %
    \item[$\bullet$] Now, Equation (\ref{eq5.34}) gives :
$$
\Phi^{n+1}(t)-\Phi^n(t)=\int_0^t\big(\sqrt{2\psi^n(s)}-\sqrt{2\psi^{n-1}(s)}\big)ds=2\int_0^t\frac{\psi^n(s)-\psi^{n-1}(s)}
{\sqrt{2\psi^n(s)}+\sqrt{2\psi^{n-1}(s)}}ds \;.
$$
In order to get rid of the denominator of the right hand side, we
use Equation (\ref{eq5.35}) which we recall is $
\dot{\psi}^{n+1}=-6U^n\psi^n-m^2\Phi^n\sqrt{2\psi^n}-\rho^2 $. Since
the sequence $(X^n)$ is uniformly bounded, on $[0,T_0]$, we have $
|-6U^n\psi^n-m^2\Phi^n\sqrt{2\psi^n}-\rho^2|\leq C_1\; $
 and thus,
$$
\frac{d\psi^{n+1}}{dt}\geq  -C_1 \;.
$$
The last inequality implies that $\psi^{n+1}(t)\geq \psi_0-C_1t \;.
$
Since $\psi_0>0$, we assume that $T_0$ is sufficiently small so that
: $0\leq C_1t\leq \frac{\psi_0}{2}, \, 0\leq t\leq T_0$  and obtain
that $\psi^{n+1}(t)\geq \frac{\psi_0}{2}$ and $
\frac{1}{\sqrt{2\psi^{n+1}(t)}}\leq\frac{1}{\sqrt{\psi_0}} $. From
there, we deduce that there exists a constant $C>0$ such that:
\begin{equation}\label{eq6.6}
\forall t\in[0,T_0], \quad  |\Phi^{n+1}(t)-\Phi^n(t)|^2\leq
C(C_0)\int_0^t|\psi^n(s)-\psi^{n-1}(s)|^2ds \;.
\end{equation}
    \item[$\bullet$] From Equation (\ref{eq5.35}) we have:
$$
\psi^{n+1}(t)-\psi^n(t)=
\int_0^t\Big[6\big(U^n(s)\psi^n(s)-U^{n-1}(s)\psi^{n-1}(s)\big)
+m^2\big(\Phi^n(s)\sqrt{2\psi^n(s)}-\Phi^{n-1}(s)\sqrt{\psi^{n-1}(s)}\big)\Big]ds
\;.
$$
Since $ \frac{1}{\sqrt{2\psi^{n+1}(t)}}\leq\frac{1}{\sqrt{\psi_0}}
$, we deduce that there exists a constant $C>0$ such that :
\bel{eq6.7}
|\psi^{n+1}(t)-\psi^n(t)|^2\leq
C\int_0^t\Big(|U^n(s)-U^{n-1}(s)|^2+|\Phi^n(s)-\Phi^{n-1}(s)|^2+|\psi^n(s)-\psi^{n-1}(s)|^2\Big)ds
 \;. \ee
    \item[$\bullet$] Finally, from Equation (\ref{eq5.36})  we get:
\begin{eqnarray}\label{aaaa}
\lefteqn{\frac{\partial (f^{n+1}-f^n)}{\partial
t}+\Big(E^nZ^n\int_{\mathbb{R}^3}f^n(t,\overline{v})d\overline{v}\Big)\sum\limits_{i=1}^3\frac{\partial(
f^{n+1}-f^n)}{\partial
u^i}}\nonumber\\
&&=\Big(E^nZ^n\int_{\mathbb{R}^3}f^n(t,\overline{v})d\overline{v}-E^{n-1}Z^{n-1}\int_{\mathbb{R}^3}f^{n-1}(t,\overline{v})
d\overline{v}\Big) \sum\limits_{i=1}^3\frac{\partial f^{n}}{\partial
u^i}
+\frac{1}{u^0_n}\tilde{Q}_n-\frac{1}{u^0_{n-1}}\tilde{Q}_{n-1}\;.\nonumber\\
\end{eqnarray}
But since the sequence $(X^n)$ is  bounded, we have:
 \begin{eqnarray*}
\lefteqn{
  \Big|E^n(t)Z^n(t)\int_{\mathbb{R}^3}f^n(t,\overline{v})d\overline{v}-E^{n-1}(t)Z^{n-1}(t)\int_{\mathbb{R}^3}f^{n-1}(t,\overline{v})d\overline{v}\Big|}
\\
 &=&\Big|(E^n(t)-E^{n-1}(t))Z^n(t)\int_{\mathbb{R}^3}f^n(t,\overline{v})d\overline{v}+(Z^n(t)-Z^{n-1}(t))E^{n-1}\int_{\mathbb{R}^3}f^n(t,\overline{v})d\overline{v}\\
&&+E^{n-1}(t)Z^{n-1}(t)\int_{\mathbb{R}^3}(f^n-f^{n-1})(t,\overline{v})d\overline{v}\Big|\\
&\leq&
C\Big(|E^n(t)-E^{n-1}(t)|+|Z^n(t)-Z^{n-1}(t)|+\|f^n(t,\cdot)-f^{n-1}(t,\cdot)\|_{H^2_d(\R^3)}\Big)
\;.
\end{eqnarray*}
Since $f^n\in H^3_{d}(\R^3)$ the second term in the r.h.s. of
Equation (\ref{aaaa})  is an element of $H^2_{d}(\R^3)$
and then,
\begin{eqnarray}\label{N1a}
\lefteqn{
\Big\|\Big(E^n(t)Z^n(t)\int_{\mathbb{R}^3}f^n(t,\overline{v})d\overline{v}-E^{n-1}(t)Z^{n-1}(t)
\int_{\mathbb{R}^3}f^{n-1}(t,\overline{v})d\overline{v}\Big)
\sum\limits_{i=1}^3\frac{\partial f^{n}}{\partial
u^i}(t,\cdot)\Big\|_{H^2_d(\R^3)}}\nonumber\\
 &\leq& C \Big(|E^n(t)-E^{n-1}(t)|+|Z^n(t)-Z^{n-1}(t)|+\|f^n(t,\cdot)-f^{n-1}(t,\cdot)\|_{H^2_d(\R^3)}\Big)
 \;.
\end{eqnarray}

Further (see inequality (\ref{N2}) of Lemma \ref{lll1}),
$$
\Big\|\frac{1}{u_n^0}\tilde{Q}_n-\frac{1}{u_{n-1}^0}\tilde{Q}_{n-1}\Big\|_{H^2_d(\R^3)}\leq
C\left(|E^n-E^{n-1}|+\|f^n-f^{n-1}\|_{H^2_d(\R^3)}\right).
$$
Equation (\ref{aaaa}) has a form to which Corollary \ref{c1}
applies. Thus,  for $k=2$ and $u= f^{n+1}-f^{n}$ in Inequality
(\ref{eq4.22}) and using the last two estimates, we have:
\begin{eqnarray}\label{eq6.9}
\lefteqn{ \|f^{n+1}(s,\cdot)-f^n(s,\cdot)\|^2_{H^2_d(\R^3)}}\\
&&\leq
C\int_0^t\Big(|E^n(s)-E^{n-1}(s)|^2+|Z^n(s)-Z^{n-1}(s)|^2+\|f^n(s,\cdot)-f^{n-1}(s,\cdot)\|^2_{H^2_d(\R^3)}\Big)ds
\;. \nonumber
\end{eqnarray}
\end{enumerate}
Consider the space
$$
\Sigma_{T_0}:=\Big(\mathcal{C}\big([0,T_0]),\R\big)\Big)^6\times\mathcal{C}\big([0,T_0],H^2_d(\R^3)\big)
\;.
$$
 Endowed with the norm
$$
|\|X|\|:=\sum\limits_{i=1}^6\sup\limits_{0\leq t\leq
T_0}|X_i(t)|+\sup\limits_{0\leq t\leq T_0}\|X_7(t)\|_{H^2_d(\R^3)}
$$
where $X=(X_i)_{1\leq i\leq 7}\in\Sigma_{T_0}$, \;$\Sigma_{T_0}$ is
a Banach space.
We want to show that there exists a constant $0<\alpha<1$ which
depends only upon $C_0$ and $T_0$  such that $|\|X^{n+1}-X^n|\| \leq
\alpha |\|X^{n}-X^{n-1}|\| $ if $T_0$ is small enough.
Summing up inequalities (\ref{eq6.2})-(\ref{eq6.7}) and
(\ref{eq6.9}) gives:
\begin{eqnarray}\label{s3.46''}
  \lefteqn{|E^{n+1}(t)-E^{n}(t)|^2+|U^{n+1}(t)-U^{n}(t)|^2+|W^{n+1}(t)-W^{n}(t)|^2+|Z^{n+1}(t)-Z^{n}(t)|^2}\nonumber\\
  &&+|\Phi^{n+1}(t)-\Phi^{n}(t)|^2 +|\psi^{n+1}(t)-\psi^{n}(t)|^2
+\|f^{n+1}(t,\cdot)-f^{n}(t,\cdot)\|^2_{H^2_d(\R^3)}\nonumber\\
   &&\leq C \int_0^t\Big(
   |E^{n}(s)-E^{n-1}(s)|^2+|U^{n}(s)-U^{n-1}(s)|^2+|W^{n}(s)-W^{n-1}(s)|^2\qquad\qquad\\
&&\;+|Z^{n}(s)-Z^{n-1}(s)|^2+|\Phi^{n}(s)-\Phi^{n-1}(s)|^2
+|\psi^{n}(s)-\psi^{n-1}(s)|^2\nonumber\\
&&\;+\|f^{n}(s,\cdot)-f^{n-1}(s,\cdot)\|^2_{H^2_d(\R^3)}\Big)ds\nonumber
\;.
\end{eqnarray}
This last inequality implies that
\begin{equation}\label{eq6.10}
|\|X^{n+1}-X^n|\|\leq \sqrt{C(C_0)T_0}|\|X^{n}-X^{n-1}|\| \;.
\end{equation}
Now, we assume that $T_0$ is small enough such that  $C(C_0)T_0<1$,
an obtain from $(\ref{eq6.10})$ that the sequence $(X^n)$ is a
Cauchy sequence in the Banach space $\Sigma_{T_0}$.
\item Next, we show that $(\frac{dY^n}{dt})$ is a Cauchy sequence in
$\big(C^0([0,T_0];\R)\big)^6$.\\
Since the sequence $(X^n)$ is bounded, from equations
(\ref{eq5.30})-(\ref{eq5.35}) we deduce that there exists a constant
$C >0$ which only depends upon $C_0$ and $T$ such that:
\begin{equation}\label{eq6.11}
    \Big\|\Big|\frac{dY^{n+1}}{dt}-\frac{dY^n}{dt}\Big\|\Big|\leq
C|\|Y^n-Y^{n-1}|\| \;.
\end{equation}

Note that  $|\|Y^n-Y^{n-1}|\|\leq |\|X^n-X^{n-1}|\|$, thus
inequality (\ref{eq6.10}) shows that:
\begin{equation}
\label{eq6.12}
    \Big\|\Big|\frac{dY^{n+1}}{dt}-\frac{dY^n}{dt}\Big\|\Big|\leq
\big(C_2\sqrt{C(C_0)T_0}\big)^{n-1}\||X^1-X^0|\| \;.
\end{equation}

As we did before, the constant $C_2\sqrt{C(C_0)T_0}$ can be assumed
to be less that one even if it means to shrink $T_0$ again. This
shows that $(\frac{dY^n}{dt})$ is a Cauchy sequence and thus,
$(Y^n)$ is a Cauchy sequence in  $\big(C^1([0,T_0],\R)\big)^6$.
\item Let us now prove that $(\frac{\partial f^n}{\partial t})$
is a Cauchy sequence in $C^0\big([0,T_0];H^1_d(\R^3)\big)$.
From Equation (\ref{aaaa})  we have:
\begin{eqnarray*}
 \Big\|\frac{\partial (f^{n+1}-f^n)}{\partial
t}\Big\|_{H^1_d(\R^3)}
  &\leq&C\Big\{\|f^{n+1}-f^n\|_{H^2_d(\R^3)}+ \|f^{n}-f^{n-1}\|_{H^1_d(\R^3)}\\
  &&
  \qquad \qquad \qquad \qquad \quad\;\;
  +\;|E^n-E^{n-1}|+|Z^n-Z^{n-1}|\Big\} \\
   &\leq& C\big\{|\|X^{n+1}-X^n|\|+|\|X^{n}-X^{n-1}|\|\big\}
   \;.
\end{eqnarray*}
Thus,
\begin{eqnarray}
\Big\| \frac{\partial (f^{n+1}-f^n)}{\partial
t}\Big\|_{C^0([0,T_0];H^1_d(\R^3))}&\leq&
C\big\{|\|X^{n+1}-X^n|\|+|\|X^{n}-X^{n-1}|\|\big\} \\
&\leq &\;C \alpha^{n}|\|X^{1}-X^{0}|\|.
\end{eqnarray}
which shows $(\frac{\partial f^n}{\partial t})$ is a Cauchy sequence
in $C^0\big([0,T_0];H^1_d(\R^3)\big)$ and ends the proof.
\end{enumerate}
\qed
\begin{remark}
Note that $(f^n)$ is a Cauchy sequence in
$C^0\big([0,T_0];H^2_d(\R^3)\big)\cap
C^1\big([0,T_0];H^1_d(\R^3)\big).$
\end{remark}
\begin{remark} \label{rq}
We notice for later use that the sequence $\left(\frac{\partial
f^n}{\partial t}\right)$ is uniformly bounded in $H^2_d(\R^3).$ In
fact we have (recall that $(X^n)$ is uniformly bounded in
$\big(\mathcal{C}^1_b([0,T_0])\big)^6\times\mathcal{C}^0\big([0,T_0];H^3_d(\R^3)\big)$):
\begin{eqnarray*}
  \Big\|\frac{\partial f^n}{\partial
t}(t,\cdot)\Big\|_{H^2_d(\R^3)} &=&
\Big\|\Big(E^{n-1}Z^{n-1}\int_{\mathbb{R}^3}f^{n-1}(t,\overline{v})d\overline{v}\Big)\sum\limits_{i=1}^3\frac{\partial
f^{n}}{\partial
u^i}+\frac{1}{u^0_{n-1}}\tilde{Q}_{n-1}\Big\|_{H^2_d(\R^3)} \\
   &\leq&\big(|E^{n-1}||Z^{n-1}|\|f^{n-1}\|_{H^3_d(\R^3)}\big)\|f^{n}\|_{H^3_d(\R^3)}+\|f^{n-1}\|^2_{H^2_d(\R^3)}\\
  &\leq&C(C_0)\;.
\end{eqnarray*}
\end{remark}
We are now ready to prove existence and uniqueness of the solution
of the Einstein-Maxwell-Boltzmann-massive scalar field system with
data  described in (\ref{eq5.19}). Before stating the main theorem
of this paper let us recall some classical facts about Sobolev
spaces $H^{(s)}$ with real $s$ (see \cite{HR}).
\subsection{Local in time solutions}
\begin{definition}
Let $n$ be a positive integer, $s$ a real number. Denote by
$\mathcal{S}'(\R^n)$ the set of all temperate distributions on
$\R^n$. We say that $u \in H^{(s)}(\R^n)$ if its Fourier transform
$\hat{u}$  is a measurable function such that
$(1+|\xi|^2)^{s/2}\hat{u}(\xi)$ is square integrable. If $u \in
H^{(s)}(\R^n)$ we define the norm
$$
\|u\|_{(s)}:= \left(
\frac{1}{(2\pi)^n}\int_{\R^n}(1+|\xi|^2)^{s}|\hat{u}(\xi)|^2d\,\xi
\right)^{1/2} \;.
$$
\end{definition}
\begin{remark} \label{rqq}
$H^{(s)}(\R^n)$ is a Hilbert space and we have
$H^{(s)}(\R^n) \hookrightarrow H^{(r)}(\R^n)$ for $s\leq r$.
Further, $H^{(0)}(\R^n)= L^{2}(\R^n)$  and $H^{(s)}(\R^n)\equiv
H^{s}(\R^n)$ when $s$ is an integer.
\end{remark}
We recall now an important inequality of functional analysis called
interpolation inequality.
\begin{lemma}
Let $s_1 < s_2 < s_3$ be real numbers and assume that $u \in
H^{(s_3)}(\R^n)$. Then  we have the following inequality
\bel{IE} \|u\|_{(s_2)}\leq
\|u\|_{(s_1)}^{\frac{s_3-s_2}{s_3-s_1}}\times\|u\|_{(s_3)}^{\frac{s_2-s_1}{s_3-s_1}}
\;. \ee
\end{lemma}
Let us state now the main theorem of the paper.
\begin{theorem} \label{MT} Let $r$ and $d$ be two positive real numbers such
that $d>5/2$. Assume that the Cauchy data
$(E_0,\,U_0,\,W_0,\,Z_0,\,\Phi_0,\,\psi_0,\,f_0)$ of the system of
partial differential equations (\ref{eq5.12}) - (\ref{eq5.18})
satisfy the Hamiltonian constraint (\ref{eq5.7}) and are such that
\begin{equation}
  \label{data}
  E_0>0, \quad W_0<0, \quad \Phi_0>0, \quad \psi_0>0, \quad f_0>0 \quad \mbox{and} \quad
  f_0\in H^3_{d,r}(\R^3)
  \;;
\end{equation}
then, there exists a positive real number $T_0$ which only depends
on $d$ and the size of the data and a unique solution
$(E,U,W,Z,\Phi,\psi,f)$ of the Cauchy problem (\ref{eq5.12}) -
(\ref{eq5.19})  defined on $[0,T_0]$  such that
$$(E,U,W,Z,\Phi,\psi,f)\in\big(\mathcal{C}^1_b([0,T_0])\big)^6\times\mathcal{C}^1\big([0,T_0]\times\R^3\big).$$
Furthermore,
\bel{cont} f\in  C^0\big([0,T_0];H^3_d(\R^3)\big)
\;. \ee
Consequently, the coupled system Einstein-Maxwell-Boltzmann-massive
scalar field has a unique local (in time) solution.
\end{theorem}
\begin{remark}
 Property (\ref{cont}) will be the key property when deriving global
 solutions for small data since we will use the continuity argument.
\end{remark}
\pr
1)\,- Existence:  As the first step towards existence of solution we
prove that the
        sequence $(X^n)$ converges in the space
        $\big(\mathcal{C}^1_b([0,T_0])\big)^6\times\mathcal{C}^1\big([0,T_0]\times\R^3\big).$
From Proposition \ref{p6.1} we know that the sequence $(Y^n)_n$ is a
Cauchy sequence in the Banach space $(C^1([0,T_0];\R))^6$ thus there
exists a set of functions $Y= (E,U,W,Z,\Phi,\psi)$ such that $(Y^n)$
converges towards $Y$ in $\big(\mathcal{C}^1_b([0,T_0])\big)^6$.
Secondly, Proposition \ref{p6.1} also tells us that
$(f^n)_n$ is a Cauchy sequence in the Banach space
$C^0\big([0,T_0];H^2_d(\R^3)\big)$. It follows that there exists a
function $f$ such that $(f^n)_n$ converges to $f$ in the space
$C^0\big([0,T_0];H^2_d(\R^3)\big)$.
 But the space
 $C^0\big([0,T_0];H^2_d(\R^3)\big)$ embeds continuously in $
C^0\big([0,T_0];H^2(\R^3)\big)$, therefore $(f^n)$ is a  Cauchy
sequence in $C^0\big([0,T_0];H^2(\R^3)\big).$ 
Now from interpolation inequality (\ref{IE}), for any real number
$2<s<3$ we have:
 \bel{a}
\|f^n(t,\cdot)-f^p(t,\cdot)\|_{(s)}\leq\|f^n(t,\cdot)-f^p(t,\cdot)\|^{s-2}_{H^3(\R^3)}
\|f^n(t,\cdot)-f^p(t,\cdot)\|^{3-s}_{H^2(\R^3)} \;. \ee
Since $\big(f^n(t, \cdot)\big)_n$ is a uniformly bounded in
$H^3_d(\R^3)$ and then in $H^3(\R^3)$, inequality $(\ref{a})$ shows
that $\big(f^n\big)_n$ is a Cauchy sequence in
$C^0\big([0,T_0];H^{(s)}(\R^3)\big)$ for any $2<s<3$.

Similarly,  since the sequence $\left(\frac{\partial f^n}{\partial
t}\right)$ is uniformly bounded (see Remark \ref{rq})  and
 is a Cauchy sequence
in $C^0\big([0,T_0];H^1_d(\R^3)\big)$, the interpolation inequality
shows that it is a Cauchy sequence in
$C^0\big([0,T_0];H^{(s)}(\R^3)\big)$ for any $1<s<2$. We then obtain
that
\bel{acs}
  (f^n) \; \mbox{is a Cauchy sequence in}\;
 C^0\big([0,T_0];H^{(s+1)}(\R^3)\big)\cap
C^1\big([0,T_0];H^{(s)}(\R^3)\big); \quad  1<s<2 . \ee
 Now, from Sobolev embedding inequality we know
that
\bel{rta} C^0\big([0,T_0];H^{(s+1)}(\R^3)\big)\cap
C^1\big([0,T_0];H^{(s)}(\R^3)\big)\hookrightarrow
C^1_b([0,T_0]\times\R^3)\quad \mbox{for any }\quad s>\frac{3}{2} \;.
\ee
Therefore choosing a particular $s$ in (\ref{acs}) such that
$\frac32 <s<2$ shows that $(f^n)$ is a  Cauchy in
$C^1_b([0,T_0]\times\R^3)$ and thus converges towards a function
$\tilde{f}$  in $C^1_b([0,T_0]\times\R^3)$ and the embedding
$C^0([0,T_0];H^2_d(\R^3))\hookrightarrow C^0([0,T_0]\times\R^3)$
shows that $f=\tilde{f}$. This shows that the collection of
functions $X=(E,U,W,Z,\Phi,\psi,f)$ is the limit of the sequence
$(X^n)$ in the space
$\big(\mathcal{C}^1_b([0,T_0])\big)^6\times\mathcal{C}^1\big([0,T_0]\times\R^3\big)$.

As the second step towards existence, we now prove that $X$ is
indeed a solution of (\ref{eq5.12}) - (\ref{eq5.18}).
Since $(Y^n)_n$  converges towards $Y= (E,U,W,Z,\Phi,\psi)$ in
$(C^1([0,T_0];\R))^6$ taking the limit pointwise in Equations
(\ref{eq5.30}), (\ref{eq5.32})-(\ref{eq5.35}) shows that $E,\,U,
\,W,\, Z,\,\Phi$ and $\psi$ satisfy (\ref{eq5.12}), (\ref{eq5.14}) -
(\ref{eq5.17}).
It remains to show that integrals
$\displaystyle\int_{\R^3}f^n(t,\overline{v})d\overline{v}$  and
$\displaystyle\int_{\R^3}\frac{|v|^2}{v_n^0}f^n(t,\overline{v})d\overline{v}$
converge respectively to
$ \displaystyle\int_{\R^3}f(t,\overline{v})d\overline{v}$ and $
\displaystyle\int_{\R^3}\frac{|v|^2}{v^0}f(t,\overline{v})d\overline{v}$,
$\forall t\in [0,T_0]$ as $n$ goes to infinity. We notice that these
last two integrals are convergent since $\forall t\in [0,T_0], \;
f(t,\cdot)\in H^2_d(\R^3)$. We have
\begin{eqnarray*}
   \Bigg|\int_{\R^3}f^n(t,\overline{v})d\overline{v}-\int_{\R^3}f(t,\overline{v})d\overline{v}\Bigg|&\leq&
\int_{\R^3}|f^n(t,\overline{v})-f(t,\overline{v})|d\overline{v}  \\
   &\leq & \|f^n(t,\cdot)-f(t,\cdot)\|_{H^2_d(\R^3)}
   \;.
\end{eqnarray*}
This shows that
$\displaystyle\int_{\R^3}f^n(t,\overline{v})d\overline{v}\longrightarrow\int_{\R^3}f(t,\overline{v})d\overline{v}$,
$\forall t\in [0,T_0]$.
 Similarly, from (\ref{N1}) we have
\begin{eqnarray*}
   \Bigg|\int_{\R^3}\frac{|v|^2}{v_n^0}f^n(t,\overline{v})d\overline{v}-
\int_{\R^3}\frac{|v|^2}{v^0}f(t,\overline{v})d\overline{v}\Bigg|
   &\leq & C\big(|E^n-E|+ \|f^n(t,\cdot)-f(t,\cdot)\|_{H^2_d(\R^3)}\big)\longrightarrow 0
   \;.
\end{eqnarray*}
Finally let us prove that
$\frac{1}{u_n^0}\tilde{Q}(f^n,f^n)\longrightarrow
\frac{1}{u^0}\tilde{Q}(f,f)$. As we did in the proof of (\ref{N2}),
we have
\begin{eqnarray*}
   \big\|\frac{1}{u_n^0}\tilde{Q}(f^n,f^n)-\frac{1}{u^0}\tilde{Q}(f,f)\big\|_{H^2_d(\R^3)}&\leq & C(T_0)
   \left(\|f^n\|_{H^2_d(\R^3)}+\|f\|_{H^2_d(\R^3)}\right)\big(|E^n-E|+\|f^n-f\|_{H^2_d(\R^3)}\big) \\
   &\leq & C(T_0)\left(C_0+\|f\|_{H^2_d(\R^3)}\right)\big(|E^n-E|+\|f^n-f\|_{H^2_d(\R^3)}\big)
   \;.
\end{eqnarray*}
Thus, $\frac{1}{u_n^0}\tilde{Q}_n$ converges towards
$\frac{1}{u^0}\tilde{Q}(f,f)$ in the space $H^2_d(\R^3)$ and since
$H^2_d(\R^3)\hookrightarrow\mathcal{C}^0_b(\R^3)$ it follows that
$\frac{1}{u_n^0}\tilde{Q}_n$ converges towards
$\frac{1}{u^0}\tilde{Q}(f,f)$ in $\mathcal{C}^0_b(\R^3)$.

We have thus proved that the limit can also be taken pointwise in
the remaining Equations (\ref{eq5.31}) and (\ref{eq5.36}) to obtain
that $X=(E,U,W,Z,\Phi,\psi,f)$ also satisfies \ref{eq5.13} and
(\ref{eq5.18}). Therefore $X$ is a local solution of (\ref{eq5.12})-
(\ref{eq5.18}).

2)-Uniqueness:
 Suppose that there exists two sets of functions
$X_i=(E_i,U_i,W_i,Z_i,\Phi_i,\psi_i,f_i),\; i=1, 2$ which solve the
system (\ref{eq5.12})-(\ref{eq5.18}) with the same Cauchy data. We
proceed exactly as we did in the proof of Proposition \ref{p6.1} and
obtain an estimate of the form (\ref{s3.46''}) where
 $X^{n+1}-X^{n}$  and $X^{n}-X^{n-1}$ are replaced by $X_2-X_1$. Applying
Gronwall's inequality to the obtained estimate proves that
$X_1=X_2.$

3)-Now, we prove that  $\forall t\in[0,T], \;
 f(t,\cdot)\in H_d^3(\R^3)$.
Recall, the sequence $\left(f^n(t,\cdot)\right)$ is uniformly
bounded in the Hilbert space $H^3_d(\R^3)$. But it is well known
that any bounded sequence in a Hilbert space has a weakly convergent
subsequence (see for example \cite{SK}, Theorem 5.4.2 page 151).
Therefore there exists a subsequence $\left(f^{n_p}(t,\cdot)\right)$
and a function $g(t, \cdot)\in H^3_d(\R^3)$ such that
$\left(f^{n_p}(t,\cdot)\right)$ converges to $g(t, \cdot)$ in
$H^3_d(\R^3)$ endowed with its weak topology which continuously
embeds in  $H^2_d(\R^3)$ endowed with its weak topology. Moreover,
recall again, $(f^n(t,\cdot))_n$ converges to $f(t,\cdot)$ in the
space $H^2_d(\R^3)$ thus this convergence also holds in
$H^2_d(\R^3)$ endowed with its weak topology. Since the weak
topology is Hausdorff, we thus conclude that $f(t, \cdot)=g(t,
\cdot)$ and then $f(t, \cdot)\in H^3_d(\R^3)$. Note that
\bel{ze} \|f(t,\cdot)\|_{H^3_d(\R^3)}\leq \liminf
\|f^{n_p}(t,\cdot)\|_{H^3_d(\R^3)} \leq C,\;\; \mbox{uniformly in}\;
t \;.\ee

4)- Let us show that (\ref{cont}) holds.
 We proceed as in \cite{HR} by first proving the weak continuity and then the strong continuity.
  The following Lemma will be needed, its
 proof is given at the end of the paper.
 \begin{lemma}
 \label{ldd}
 The space $\mathcal{C}^{\infty}_c(\R^3)$ of Compactly supported smooth functions defined on
  $\R^3$ is dense in the space $H^3_d(\R^3)$.
 \end{lemma}
 \noindent
\textbf{Weak continuity:}
 First, let us prove that the
solution is weakly continuous. Let $F$ be in the dual of
$H^3_d(\R^3)$. Then, due to the Riesz representation theorem there
exists   $\varphi_F\in H^3_d(\R^3)$ such that for all $h\in
H^3_d(\R^3)$,
$$
F(h)=\langle h,\varphi\rangle_{H^3_d(\R^3)}=\sum_{|\alpha|\leq
3}\int_{\R^3}(1+|\bar{v}|)^{d+|\alpha|}\partial^{\alpha}h(\bar{v})\cdot
(1+|\bar{v}|)^{d+|\alpha|}\partial^{\alpha}\varphi_F(\bar{v})d\bar{v}.
$$
We then have,
$$
F(f^n(t,\cdot))-F(f(t,\cdot))=\sum_{|\alpha|\leq
3}\int_{\R^3}(1+|\bar{v}|)^{d+|\alpha|}\partial^{\alpha}f^n\cdot(1+|\bar{v}|)^{d+|\alpha|}\partial^{\alpha}\varphi_F
d\bar{v} - \sum_{|\alpha|\leq
3}\int_{\R^3}(1+|\bar{v}|)^{d+|\alpha|}\partial^{\alpha}f\cdot(1+|\bar{v}|)^{d+|\alpha|}\partial^{\alpha}\varphi_F
d\bar{v} \;.
$$
Consequently, if $(\varphi_j)$ is a sequence of compactly supported
 smooth functions converging to $\varphi_F$ in $H^3_d(\R^3)$, we
obtain:
\begin{eqnarray*}
 F(f^n(t,\cdot))-F(f(t,\cdot))  &=&
\sum_{|\alpha|\leq
3}\int_{\R^3}(1+|\bar{v}|)^{d+|\alpha|}\partial^{\alpha}(f^n-f)(t,\,\bar{v})\cdot(1+|\bar{v}|)^{d+|\alpha|}\partial^{\alpha}
(\varphi_F-\varphi_j)(\bar{v})d\bar{v} \\
&&+\; \sum_{|\alpha|\leq
3}\int_{\R^3}(1+|\bar{v}|)^{d+|\alpha|}\partial^{\alpha}(f^n-f)(t,\,\bar{v})\cdot(1+|\bar{v}|)^{d+|\alpha|}\partial^{\alpha}
\varphi_j(\bar{v})d\bar{v} \;.
\end{eqnarray*}
We then deduce since, the sequence $(f^n)$ is bounded in
$H^3_d(\R^3)$ that:
$$
\big|F(f^n(t,\cdot))-F(f(t,\cdot))\big|\leq
C\|\varphi_F-\varphi_j\|_{H^3_d(\R^3)}+\bigg| \sum_{|\alpha|\leq
3}\int_{\R^3}(1+|\bar{v}|)^{d+|\alpha|}\partial^{\alpha}(f^n-f)(t,\,\bar{v})\cdot(1+|\bar{v}|)^{d+|\alpha|}\partial^{\alpha}
\varphi_j(\bar{v})d\bar{v}\bigg|
$$
Letting $j$ be large enough that the first term on the right-hind
side is less than or equal to $\frac{\varepsilon}{2}$, and then
choosing $n$ large enough, depending on $j$, so that the second term
is less than $\frac{\varepsilon}{2}$, we conclude that the
right-hand side is less than $\varepsilon$. We conclude that
$F(f^n(t,\cdot))$ converges uniformly to $F(f(t,\cdot))$  which
proves that the solution $f$ is weakly continuous.
\\
\noindent \textbf{Strong continuity:} Let $t_0\in [0, T_0)$. We want
to prove that $f : [0,T_0)\longrightarrow H_d^{3}(\R^3)$ is
continuous at $t_0$ i.e. $$ \lim\limits_{t\rightarrow
t_0}\|f(t,\cdot)-f(t_0,\cdot)\|_{H_d^{3}(\R^3)}=0\;. $$
\noindent Using the inner product $\langle\, ,\,\rangle$ on
$H_d^{3}(\R^3)$, we can write for $t\in[0,T_0)$:
\be\label{sc} \langle
f(t,\cdot)-f(t_0,\cdot),f(t,\cdot)-f(t_0,\cdot)\rangle=\langle
f(t,\cdot),f(t,\cdot)\rangle-2\langle
f(t,\cdot),f(t_0,\cdot)\rangle+\langle
f(t_0,\cdot),f(t_0,\cdot)\rangle\;. \ee
Note that the last term on the right-hand side is
$\|f(t_0,\cdot)\|^2_{H_d^{3}(\R^3)}$. Due to the weak continuity of
$f$, the limit as $t$ goes to $t_0$ of second term on the right-hand
side is
$-2\|f(t_0,\cdot)\|^2_{H_d^{3}(\R^3)}$.\\
For the first term on the right-hand side, we suppose that $t>t_0$
and use the fact that there exists $\delta>0$ such that (see
estimate (\ref{eq4.22}))
$$
e^{-\delta (t-t_0)}\|f(t,\cdot)\|^2_{H^3_d(\R^3)}\leq
\|f(t_0,\cdot)\|^2_{H^3_d(\R^3)}+ C\int_{t_0}^te^{-\delta
s}\|f(s,\cdot)\big\|^4_{H^3_d(\R^3)}ds\;;
$$
and the fact that $f(t,\cdot)$ is uniformly bounded (see
(\ref{ze}))), to have
$ \lim\limits_{t\rightarrow t_0^+} \langle
f(t,\cdot),f(t,\cdot)\rangle\leq
\|f(t_0,\cdot)\|^2_{H_d^{3}(\R^3)}\,. $
Combine these observations with (\ref{sc}), to have
$ \lim\limits_{t\rightarrow t_0^+} \langle
f(t,\cdot)-f(t_0,\cdot),f(t,\cdot)-f(t_0,\cdot)\rangle\leq 0\;, $
and conclude, since $\langle
f(t,\cdot)-f(t_0,\cdot),f(t,\cdot)-f(t_0,\cdot)\rangle\geq
0,\;\forall t\in[0,T_0)$ that
$ \; \lim\limits_{t\rightarrow t_0^+} \langle
f(t,\cdot)-f(t_0,\cdot),f(t,\cdot)-f(t_0,\cdot)\rangle= 0 \;; $ i.e.
$f : [0,T_0)\longrightarrow H_d^{3}(\R^3)$ is right continuous on
$[0,T_0)$. By time reversal one obtains left continuity and thus
continuity of $f$.

 5)- Finally, since the Einstein-Maxwell-Boltzmann-massive scalar
field system is equivalent to the system of first order partial
differential equations (\ref{eq5.12})-(\ref{eq5.18}), we have
therefore proved that the Einstein-Maxwell-Boltzmann-massive scalar
field equations have a unique local (in time) solution.
\qed

\section{Global existence of solutions to the Einstein-Maxwell-Boltzmann system with massive scalar field}
\subsection{The method}
In this section, we prove under further assumptions on the data that
the local solution obtained in Section \ref{s6} is in fact a global
solution. We will use the well known continuation criterium (see for
example \cite{MiTa3}, Proposition 1.5 p. 365) which says that, the
breakdown of a classical solution $u$ of a system of hyperbolic
partial differential equations must involve a blow-up of either
$\sup\limits_{x}|u(t,x)|$ or $\sup\limits_{x}|\nabla_x u(t,x)|$. In
other words, if the $C^1$-norm of a solution $u(t, \cdot)$ on an
interval $[0,T)$ is uniformly bounded then this solution can be
extended beyond $T$.
Let us sketch out the method we adopt. Denote $[0,\,T_*)$, $T_*>0$,
the maximal existence time interval of the solution of system
 (\ref{eq5.12})-(\ref{eq5.18}) , with initial
data defined by $(\ref{eq5.19})$, subject to the Hamiltonian
constraint $(\ref{eq5.7})$ and satisfying (\ref{data}).
Assume by contradiction that $\ts<\infty$ (otherwise $\ts=+\infty$
and there is nothing to do). Then we will prove using a continuity
type argument that the solution $(E, U, W, Z, \Phi,
    \psi, f)$ is uniformly bounded on $[0,\,\ts)$ by a constant depending
    only on the initial data, $\ts ,\, m,\, r$ and $\Lambda$. It will then
follow by the continuation criterium that this solution  can be
extended to a larger time interval $[0, T')$ thus contradicting the
maximality of $\ts$. This will imply that $\ts = +\infty$ and the
solution is global. Before doing this, let us give some useful
estimates on the obtained local solution.
%
%

\subsection{A priori estimates and global solution}
\begin{lemma} \label{l2}
In addition to hypotheses of Theorem \ref{MT} assume $U_0>0$ and
that the cosmological constant satisfies $\Lambda > -4\pi
m^2\Phi_0.$ Then the solution $(E,U,W,Z,\Phi,\psi)$ defined on $[0,
\ts)$ of system (\ref{eq5.12})-(\ref{eq5.17}) with Cauchy data given
in (\ref{data}) satisfies the set of inequalities
\begin{equation}\label{Eq7.1}
\left\{
  \begin{array}{ll}
    \bigg(\frac{\Lambda}{3}+\frac{4\pi}{3}m^2\Phi^2(0)\bigg)^{1/2}\leq U(t)\leq U_0; \;
0\leq E(t)\leq E_0; \;\;0\leq
\frac{1}{E}\leq\frac{1}{E_0}e^{U_0T}& \hbox{} \\
  |Z(t)|\leq |Z_0|, \;\;
|\Phi|\leq \left(\frac{3U_0^2-\Lambda}{4\pi m^2}\right)^{1/2};\;
\;0\leq \psi\leq \frac{3U_0^2-\Lambda}{8\pi}; \;\;|W|\leq
|W_0|+\frac{\rho^2}{3U_0}e^{3U_0T} & \hbox{}
  \end{array}
  \;;
\right.
\end{equation}
For all $T$ such that $0<T<\ts$.
\end{lemma}
\pr
With the change of functions $(\ref{eq5.8})$, the Hamiltonian
constraint $(\ref{eq5.1})$ writes :
\begin{equation}\label{eq7.2}
3U^2-\Lambda=8\pi
E^3\int_{\R^3}v^0f(t,\overline{v})d\overline{v}+12\pi\frac{Z^2}{E^2}-8\pi
W+4\pi(2\psi+m^2\Phi^2) \;.
\end{equation}
Taking $(-2)\times(\ref{eq5.13})+(\ref{eq7.2})$ gives :
 \begin{equation}\label{eq7.3}
    \dot{U}=-4\pi\Big[E^5\int_{\R^3}\frac{(v^1)^2}{v^0}f(t,\overline{v})d\overline{v}
    +E^3\int_{\R^3}v^0f(t,\overline{v})d\overline{v}\Big]-8\pi\frac{Z^2}{E^2}
    +4\pi W-8\pi\psi
    \;.
 \end{equation}
 But since by $(\ref{eq5.20})$ we have $W\leq 0,\, $  $(\ref{eq7.3})$
 implies that $\dot{U}\leq 0$.
 So $U$ is decreasing. The Hamiltonian constraint $(\ref{eq7.2})$
 implies, since $\Phi^2$ is increasing :
  %
  %
 \begin{equation}\label{eq7.4}
    U^2\geq \frac{\Lambda}{3}+\frac{4\pi}{3}m^2\Phi_0^2
    \;.
 \end{equation}
 But by hypothesis,
 $
    \Lambda\geq-4\pi m^2\Phi_0^2
 $
 thus,
$(\ref{eq7.4})$ is equivalent to :
$$
\Big(U +\sqrt{\frac{\Lambda}{3}+\frac{4\pi}{3}m^2\Phi^2(0)}\Big)
\Big(U-\sqrt{\frac{\Lambda}{3}+\frac{4\pi}{3}m^2\Phi^2(0)}\Big)\geq
0 \;;
$$
which implies :
$$
U\leq-\sqrt{\frac{\Lambda}{3}+\frac{4\pi}{3}m^2\Phi^2(0)}\qquad
\mbox{or}\qquad
U\geq\sqrt{\frac{\Lambda}{3}+\frac{4\pi}{3}m^2\Phi^2(0)} \;.
$$
Now, from our hypotheses, $U_0>0$ and since $U$ is a continuous
function we will only have:
\begin{equation}\label{eq7.6}
    U\geq\sqrt{\frac{\Lambda}{3}+\frac{4\pi}{3}m^2\Phi^2(0)}
    \;,
\end{equation}
and then ($U$ is decreasing):
\begin{equation}\label{eq7.7}
    \sqrt{\frac{\Lambda}{3}+\frac{4\pi}{3}m^2\Phi^2(0)}\leq U\leq U_0
    \;.
\end{equation}
By $(\ref{eq7.7})$, $U$ is bounded. Now we have by $(\ref{eq5.12})$
: $\dot{E}=-UE$ since $E=\frac{1}{a}\geq 0$ and $U> 0$, this implies
$\dot{E}< 0$ so $E$ is decreasing and :
\begin{equation}\label{eq7.8}
    0\leq E\leq E_0=\frac{1}{a_0}
    \;.
\end{equation}
Now by $(\ref{eq5.15})$ we have $\dot{Z}=-3UZ$; since $U\geq 0$ is
bounded, we deduce that :
\begin{equation}\label{eq7.9}
    |Z(t)|\leq |Z_0|
    \;.
\end{equation}
Moreover, Equation \ref{eq5.12} gives after integration on $[0,t],
\, t>0$:
$$
0\leq \frac{1}{E(t)}\leq
\frac{1}{E_0}e^{U_0t}\leq\frac{1}{E_0}e^{U_0T}
\;.
$$
Further, the Hamiltonian constrain $(\ref{eq7.2})$ implies
$
8\pi\psi\leq 3U^2-\Lambda$  and  $4\pi m^2\Phi^2\leq 3U^2-\Lambda .
$
Since  $0\leq U(t)\leq U_0$  and $\Phi, \, \psi>0$, we have
$$
0< \psi \leq \frac{3U_0^2-\Lambda}{8\pi}\qquad \mbox{and} \qquad
                 0<\Phi \leq \sqrt{\frac{3U_0^2-\Lambda}{4\pi m^2}}
\;.
$$
Finally, integrating on $[0,\,t]$ Equation $(\ref{eq5.14})$ gives:
$$
W(t)=e^{-\int_0^t3U(s)ds}\left\{W(0)-\rho^2\int_0^te^{\int_0^s3U(\tau)d\tau}ds\right\}\;.
$$
Once more we use  $0\leq U\leq U_0$ and get
$ |W(t)|\leq |W(0)|+\rho^2\int_0^te^{3U_0s}ds\leq
|W_0|+\frac{\rho^2}{3U_0}e^{U_0T} \;; $
and the proof is complete. \qed

We are now ready to state and prove existence of a unique global
solution of the Einstein-Maxwell-Boltzmann equations with massive
scalar field and cosmological constant.

\begin{theorem}
\label{MTT} The hypotheses are those of Lemma \ref{l2}. If the real
number  $r$ (the norm of $f_0$ in the space $H^3_d(\R^3)$) is small
enough then the local solution obtained in Theorem \ref{MT} is
global in time. Consequently, the Einstein-Maxwell-Boltzmann
equations with massive scalar field and cosmological constant have a
unique global (in time) solution for small data.
\end{theorem}
\pr
As we mentioned earlier it will suffice to prove that the $C^1$-norm
of the set of functions $X= (E, U, W, Z, \Phi, \psi, f)$ is
uniformly bounded on $[0,\,\ts).$ Since from Lemma \ref{l2}, the
function $Y= (E, U, W, Z, \Phi, \psi)$ is bounded it will suffice to
prove that if $r$ is sufficiently small then there exists a constant
$M>0$ such that
$$
\sup_{t\in [0,\ts)}\|e^{-\delta_1t/2}f(t,\cdot)\|_{H^3_d(\R^3)}\leq
M \;.
$$
We use a continuity argument.
By hypothesis $f_0\in  H^3_{d,r}(\R^3)$ thus
$\|f_0\|_{H^3_d(\R^3)}\leq r$ and by continuity
\bel{BH}e^{-\delta_1t/2}\|f(t,\cdot)\|_{H^3_d(\R^3)}\leq 4r \ee
on a sub-interval of $[0,\ts).$  Denote by $[0, T_1]$ the largest
time interval on which (\ref{BH}) still holds and let us show that
if $r$ is sufficiently small then on the interval $[0, T_1]$
inequality (\ref{BH}) implies the same inequality with the constant
$4$ replaced by 2. It will then follow by continuity that there
exists a real number $\epsilon>0$ such that (\ref{BH}) still holds
on $[0, T_1+\epsilon]$ which contradicts the maximality of $T_1$ and
then $T_1=\ts$.

Recall $f$ is a $C^1-$solution of the hyperbolic Equation
(\ref{eq5.18}) which satisfies hypotheses of Corollary \ref{c1} with
$n=k=3$. Therefore inequality (\ref{eq4.22}) reads instead:
$$
e^{-\delta_1t}\|f(t,\cdot)\|^2_{H^3_d(\R^3)}\leq
 \|f_0\|^2_{H^3_d(\R^3)}+C_1\int_0^t
e^{-\delta_1s}\left\|\frac{1}{u^0}\tilde{Q}(f,f)\right\|^2_{H^3_d(\R^3)}ds
\;.
$$
This last inequality implies, using (\ref{eq4.31}) and (\ref{BH})
that
$$
e^{-\delta_1t}\|f(t,\cdot)\|^2_{H^3_d(\R^3)}\leq
\|f_0\|^2_{H^2_d(\R^3)}+16 C_1r^2 \int_0^t e^{\delta_1s}\cdot
e^{-\delta_1s}\|f(s,\cdot)\|^2_{H^3_d(\R^3)}ds  \;.
$$
Using again Gronwall's inequality, one has:
$$
e^{-\delta_1t}\|f(t,\cdot)\|^2_{H^3_d(\R^3)}\leq
\|f_0\|^2_{H^3_d(\R^3)}e^{16C_1r^2
\frac{e^{\delta_1T_0}-1}{\delta_1}}\leq
 r^2e^{16C_1r^2
\frac{e^{\delta_1T_0}-1}{\delta_1}} \;.
$$
Note that $\lim\limits_{r\,\rightarrow \,0^+}e^{16C_1r^2
\frac{e^{\delta_1T_0}-1}{\delta_1}}= 1$ thus there exists a small
$r_0>0$ such that
$$
0< r \leq r_0 \quad \Longrightarrow e^{16C_1r^2
\frac{e^{\delta_1T_0}-1}{\delta_1}}\leq 4\qquad \mbox{i.e.}\quad
e^{-\delta_1t/2}\|f(t,\cdot)\|_{H^3_d(\R^3)}\leq 2r \;.
$$
Now fix a real number $r$ such that $0< r \leq r_0$ and pick $f_0$
in $H^3_{d,r}(\R^3)$ then $\sup\limits_{t\in[0, T_1]}e^
{-\delta_1t/2}\|f(t,\cdot)\|_{H^3_d(\R^3)}\leq 2r$. By continuity,
there exists a real number $\epsilon>0$ such that
$\sup\limits_{t\in[0,
T_1+\epsilon]}e^{-\delta_1t/2}\|f(t,\cdot)\|_{H^3_d(\R^3)}\leq 4r$
which contradicts the maximality of $T_1$ and then $T_1= \ts$. We
have thus proved that
\bel{ifi} %
 \sup_{t\in
[0,\ts)}\|e^{-\delta_1t/2}f(t,\cdot)\|_{H^3_d(\R^3)}\leq 4r \;.
 \ee
This proves that the $C^1$-norm of $f$ and then of  $X=(E, U, W, Z,
\Phi,\psi, f)$ is uniformly bounded on $[0,\,\ts)$ and from the
continuation criterium can be extended as a $C^1-$solution beyond
$\ts$ which in turn contradicts the maximality of $\ts$ and thus
$\ts=+\infty$ i.e. the solution is global.
\qed
\begin{remark}
We have proved a global in time existence and uniqueness theorem for
the coupled Einstein-Maxwell-Boltzmann system which governs the
dynamics of charged particles in the presence of a massive scalar
field. In an up coming work we will study the asymptotic behavior,
the stability of the obtained solution, the geodesic completeness of
the solution, and the energy conditions in the constructed
space-space time. See \cite{NNAN2, HL2}.

\end{remark}

\appendix

\section{Energy estimate for a first order hyperbolic partial
differential equation}\label{Apend}  Let us consider the first order
PDE in $u=u(t,x)$ :
\begin{equation}\label{eq4.1}
    u_t+\sum\limits_{i=1}^na_i(x,t)u_{x_i}+b(x,t)u=f(x,t)\;\;in\;\;
    \R^n\times\R
    \;;
\end{equation}
with initial data :
\begin{equation}\label{eq4.2}
    u(x,0)=u_0(x)\;\;in\;\;
    \R^n\;.
\end{equation}
Here $b$ and $f$ are functions defined on $\R^{n+1}$ and
$a=(a_1,\cdots,a_n)$ a family of functions such that:
\begin{equation}\label{eq4.3}
   \sum\limits_{i=1}^n\big(\sup\limits_{(x,t)}|a_i(x,t)|\big)=:|a|\leq\frac{1}{\kappa}\;,
\end{equation}
where $\kappa$ is a positive constant.
\begin{remark}
We would like the point out the fact that the notations in this
subsection are independent of those of the other sections of the
paper. For examples, the letter $a$ here is used for a collection of
real valued functions and has nothing to do with the expansion
factor of the previous sections, $f$ here is the source term of the
PDE we are dealing with and must not be confused with the
distribution function. We hope that this clash of notations will not
confuse the reader.
\end{remark}

\begin{proposition} \label{p1}
Let $a=(a_1,\cdots,a_n)$ be a family of class $C^1$ functions with
bounded partial d\'erivatives with respect to $x_i$ defined on
$\R^{n+1}$ and satisfying $(\ref{eq4.3})$. Let $b$ be a bounded
function defined on  $\R^{n+1}$, and $u$ a solution of the initial
value problem $(\ref{eq4.1})-(\ref{eq4.2})$.
 Then, for every $T>0$, if
$f\in C\big([0,T);L^2(\R^n)\big)$ and  $u_0\in L^2(\R^n)$, we have:
\begin{equation}\label{eq4.4}
\int_{\R^n}e^{-\alpha t}u^2dx\leq \int_{\R^n} u^2_0
dx+\int_0^te^{-\alpha s}\|f(s,\cdot)\|^2_{L^2(\R^n)}ds \;,\qquad
t\leq T;
\end{equation}
where
 $\alpha$ is a positive constant.
\end{proposition}
\pr
 Let $\overline{t}>0$ be given. Define
$D=D_{\;\kappa,t,\overline{t}\;}\;,0\leq t<T<\overline{t}$ by :
\begin{equation}\label{eq4.6}
    D=\{(x,s)\in\R^{n+1};\;\kappa|x|<\overline{t}-s,\;0<s<t\}
    \;.
\end{equation}
Let us denote by $\Sigma_{t,\overline{t}}$,
$\Sigma_{0,\overline{t}}$ and $S_{t,\overline{t}}$ respectively, the
upper boundary, the lower boundary and the side of $D$; ie :
\begin{equation}\label{eq 7}
    \left\{
      \begin{array}{ll}
       \Sigma_{t,\overline{t}}= \{(x,t)\in\R^{n+1};\;\kappa|x|<\overline{t}-t\}& \\
        \Sigma_{0,\overline{t}}= \{(x,0)\in\R^{n+1};\;\kappa|x|<\overline{t}\} &\\
         S_{t,\overline{t}}= \{(x,s)\in\R^{n+1};\;\kappa|x|=\overline{t}-s,\;0<s<t\}&
      \end{array}
    \right.
    \;.
\end{equation}
For $\alpha>0$, we multiply Equation  \ref{eq4.1} by $2e^{-\alpha
t}u$, and obtain an equation which can be written as:
\bel{eq4.9}
    \big(e^{-\alpha t}u^2\big)_t+\sum_{i=1}^n\big(e^{-\alpha
t}a_iu^2\big)_{x_i}+e^{-\alpha
t}\big(\alpha+2b-\sum_{i=1}^na_{i,x_i}\big)u^2=2e^{-\alpha t}u\,f
\;. \ee
The first two terms of (\ref{eq4.9}) can be written as a divergence.
If we set $X = X^{\mu}\partial_{\mu}:= e^{-\alpha t}u^2 a^{\mu}
\partial_\mu,$ then
$$ div X=\sum\limits_{\mu=0}^n\big(e^{-\alpha
t}a^{\mu}u^2\big)_{x_{\mu}}  = \big(e^{-\alpha
t}u^2\big)_t+\sum\limits_{i=1}^n\big(e^{-\alpha t}a_iu^2\big)_{x_i}
\;;
$$
with $a^{\mu}=(1, a_1, \ldots, a_n)$ and $x_0=t$.
We integrate (\ref{eq4.9}) on $D$. We have
\bel{eq4.9a}
    \int_D\left(\big(e^{-\alpha t}u^2\big)_t+\sum_{i=1}^n\big(e^{-\alpha
t}a_iu^2\big)_{x_i}+e^{-\alpha
t}\big(\alpha+2b-\sum_{i=1}^na_{i,x_i}\big)u^2\right)dt\,dx=2\int_D
e^{-\alpha t}u\, f dt\,dx \;. \ee
 By the
Stokes theorem, we have :
$$
\int_{D}\Big[\big(e^{-\alpha
t}u^2\big)_t+\sum\limits_{i=1}^n\big(e^{-\alpha
t}a_iu^2\big)_{x_i}\Big]dxdt = \int_{\partial D}X\cdot\eta dS=
 \int_{\Sigma_{0,\overline{t}}}X\cdot\eta
    dS+\int_{\Sigma_{t,\overline{t}}}X\cdot\eta
    dS+\int_{S_{t,\overline{t}}}X\cdot\eta
    dS
    \;,
$$
where $\eta$ is the outward unit normal vector to $\partial D$ and
$dS$ the surface element on $\partial D$.
But,
\begin{enumerate}
\item[$\bullet$] on $S_{t,\overline{t}}, \;
\eta=(\eta_t,\eta_1,\cdots,\eta_n)=\frac{1}{\sqrt{1+\kappa^2}}\left(1,\kappa\frac{x_1}{|x|},
\cdots,\kappa\frac{x_n}{|x|} \right)$
from where we have
$$\int_{S_{t,\overline{t}}}X\cdot\eta
dS=\int_{S_{t,\overline{t}}}e^{-\alpha
t}\big(\sum\limits_{i=1}^n\eta_ia_i+\eta_t\big)u^2dS \;.
$$
\item[$\bullet$] on $\Sigma_{t,\overline{t}}$ and
$\Sigma_{0,\overline{t}}$ the outward unit normals are respectively
$\eta=(1,\cdots,0,0)$ and $\eta=(-1,\cdots,0,0)$,
and then
$$
\int_{\Sigma_{t,\overline{t}}}X\cdot\eta
dS=\int_{\Sigma_{t,\overline{t}}}e^{-\alpha t}u^2dx \qquad
\mbox{and}\qquad \int_{\Sigma_{0,\overline{t}}}X\cdot\eta
dS=-\int_{\Sigma_{0,\overline{t}}}u_0^2dx \;.
$$
\end{enumerate}

Let us observe that, by  $(\ref{eq4.3})$ we have
$\sum\limits_{i=1}^n\eta_ia_i+\eta_t\geq 0$ which leads to the
following inequality:
    \begin{equation}\label{eq4.10}
 \int_{D}\Big[\big(e^{-\alpha t}u^2\big)_t+\sum\limits_{i=1}^n\big(e^{-\alpha
t}a_iu^2\big)_{x_i}\Big]dxdt\geq\int_{\Sigma_{t,\overline{t}}}e^{-\alpha
t}u^2dx-\int_{\Sigma_{0,\overline{t}}}u_0^2dx \;.
    \end{equation}

Now since $b$ and the partial derivatives of $a_i$ with respect to
$x_i$ are bounded, we can choose $\alpha$ such that:
\begin{equation}\label{eq4.11}
\alpha+2b-\sum_{i=1}^na_{i,x_i}\geq 1\;\;\mbox{in}\;\;D \;.
\end{equation}
From now, we suppose that $\alpha$ is chosen such that
(\ref{eq4.11}) holds. We then have
\bel{eq4.11a} \int_{D}e^{-\alpha
t}\big(\alpha+2b-\sum_{i=1}^na_{i,x_i}\big)u^2 dtdx\geq
\int_{D}e^{-\alpha t}u^2 dtdx \;. \ee
Now, recall the trivial inequality
     \begin{equation}\label{eq4.12}
\int_{D}2e^{-\alpha t}ufdxdt\leq \int_{D}e^{-\alpha
t}u^2dxdt+\int_{D}e^{-\alpha t}f^2dxdt \;.
\end{equation}
Adding (\ref{eq4.10}) and (\ref{eq4.11a}) and using (\ref{eq4.9}),
(\ref{eq4.12}) give:

\begin{equation} \label{eq4.13}
    \int_{\Sigma_{t,\overline{t}}}e^{-\alpha t}u^2dx\leq
    \int_{\Sigma_{0,\overline{t}}}u_0^2dx+\int^t_0ds\int_{\Sigma_{s,\overline{t}}}e^{-\alpha
    s}f^2(s,\cdot)dx
    \;;
\end{equation}
where in (\ref{eq4.13}),  $ \Sigma_{s,\overline{t}}=
\{(x,s)\in\R^{n+1};\;\kappa|x|<\overline{t}-s\,\}$ for fixed $s$.
Now, we let $\overline{t}$ tends to $ \infty$, then
$\Sigma_{t,\overline{t}}$, $\Sigma_{s,\overline{t}}$ et
$\Sigma_{0,\overline{t}}$ tend respectively to $\{t\}\times\R^n$,
$\{s\}\times\R^n$ and $\{0\}\times\R^n$. Therefore $(\ref{eq4.13})$
gives:
\begin{equation}\label{eq4.14}
\int_{\R^n}e^{-\alpha t}u^2dx\leq
    \int_{\R^n}u_0^2dx+\int^t_0ds\int_{\R^n}e^{-\alpha
    s}f^2(s,\cdot)dx
 \;;
\end{equation}
which is the desired inequality.  \qed

\begin{remark}
For $0<t\leq T$, inequality $(\ref{eq4.14})$ reads
\begin{equation}\label{eq4.17}
e^{-\alpha t}\|u(t,\cdot)\|^2_{L^2(\R^n)}\leq
\|u_0\|^2_{L^2(\R^n)}+\int_0^t e^{-\alpha
s}\|f(s,\cdot)\big\|^2_{L^2(\R^n)}ds \;.
\end{equation}
\end{remark}
The $H^k-$version of the previous Proposition reads:

\begin{proposition} \label{ta}
Let $a=(a_1,\cdots,a_n), b$ be a family of smooth functions defined
on $\R^{n+1}$ such that the collection of functions $a$ satisfies
$(\ref{eq4.3})$. Suppose that the partial derivatives up to order
$k\in\N^{\ast}$ of $a$ and $b$ with respect to the space variables
$x$ are bounded. Assume that $u$ is a $C^1$ solution of the initial
value problem $(\ref{eq4.1})-(\ref{eq4.2})$.
 Then, for every $T>0$, if
$f\in C\big([0,T);H^k(\R^n)\big)$ and  $u_0\in H^k(\R^n)$, then for
all $t\in[0,T], \; u(t, \cdot)\in H^k(\R^n)$ and
\begin{equation}\label{eq4.18}
e^{-\delta_0 t}\|u(t,\cdot)\|^2_{H^k(\R^n)}\leq
\|u_0\|^2_{H^k(\R^n)}+ C_0\int_0^t\|e^{-\delta_0
s}f(s,\cdot)\big\|^2_{H^k(\R^n)}ds \;;
\end{equation}
where $\delta_0$ and  $C_0$ are positive constants.
\end{proposition}
\pr
Note that $(\ref{eq4.17})$ gives  $(\ref{eq4.18})$ for $k=0$.
Consider the commutator $[L,\,\partial^{\alpha}]$ defined by :
$[L,\,\partial^{\alpha}]u=L\partial^{\alpha}u-\partial^{\alpha}Lu$,
where $L$ is the  differential operator associated to
$(\ref{eq4.1})$ (here,  $\partial^{\alpha}\equiv
\partial^{\alpha}_x)$. We have, since $(\ref{eq4.1})$:
\begin{equation}\label{eq4.19}
L\partial^{\alpha}u=\partial^{\alpha}f+[L,\,\partial^{\alpha}]u \;.
\end{equation}
Applying Proposition \ref{p1} to (\ref{eq4.19}) with  $u$ replaced
by $\partial^{\alpha}u$ and $f$ replaced by
$\partial^{\alpha}f+[L,\,\partial^{\alpha}]u$ shows that there
exists $\delta>0$ such that:
\begin{eqnarray}
 e^{-\delta
t}\|(\partial^{\alpha}u)(t,\cdot)\|^2_{L^2(\R^n)}&\leq&
 \|(\partial^{\alpha}u)(0,\cdot)\|^2_{L^2(\R^n)}+ \int_0^te^{-\delta
s}\big\|(\partial^{\alpha}f+[L,\,\partial^{\alpha}]u)(s,\cdot)\big\|^2_{L^2(\R^n)}ds
\nonumber
\\
&\leq&
\|(\partial^{\alpha}u)(0,\cdot)\|^2_{L^2(\R^n)}+ 2\int_0^te^{-\delta
s}\big\|\partial^{\alpha}f\big\|^2_{L^2(\R^n)}+2\int_0^te^{-\delta
s}\|[L,\,\partial^{\alpha}]u(s,\cdot)\big\|^2_{L^2(\R^n)}ds
\;; \nonumber
\\
 \label{eq4.20}
\end{eqnarray}
where $|\alpha|\leq k$.
 On the order hand,
\begin{eqnarray}
[L,\,\partial^{\alpha}]u&=&L\partial^{\alpha}u-\partial^{\alpha}Lu
\nonumber\\
&=&
(\partial^{\alpha}u)_t+\sum\limits_{i=1}^na_i(\partial^{\alpha}u)_{x_i}+b\partial^{\alpha}u-\partial^{\alpha}
u_t-\partial^{\alpha}\big(\sum\limits_{i=1}^na_iu_{x_i}\big)-\partial^{\alpha}(bu)
\nonumber
\\
&=&
\sum\limits_{i=1}^na_i(\partial^{\alpha}u)_{x_i}+b\partial^{\alpha}u
-\sum\limits_{i=1}^n\sum\limits_{0\leq\beta\leq\alpha}\complement^{\beta}_{\alpha}\partial^{\beta}
a_i\partial^{\alpha-\beta}u_{x_i}-\sum\limits_{0\leq\beta\leq\alpha}
\complement^{\beta}_{\alpha}\partial^{\beta}b\partial^{\alpha-\beta}u
\nonumber
\\
&=&
-\sum\limits_{i=1}^n\sum\limits_{0<\beta\leq\alpha}\complement^{\beta}_{\alpha}\partial^{\beta}a_i\partial^{\alpha-\beta}u_{x_i}
   -\sum\limits_{0<\beta\leq\alpha}\complement^{\beta}_{\alpha}\partial^{\beta}b\partial^{\alpha-\beta}u\;.
   \label{eq4.21}
\end{eqnarray}
Since $\partial^{\beta}a_i$ and $\partial^{\beta}b$ are bounded, the
last equality gives the following estimate
\begin{eqnarray*}
\|[L,\,\partial^{\alpha}]u\|_{L^2(\R^n)}&\leq&
   \Big\|\sum\limits_{i=1}^n\sum\limits_{0<\beta\leq\alpha}\complement^{\beta}_{\alpha}
   \partial^{\beta}a_i\partial^{\alpha-\beta}u_{x_i}\Big\|_{L^2(\R^n)}+
   \Big\|\sum\limits_{0<\beta\leq\alpha}\complement^{\beta}_{\alpha}\partial^{\beta}
   b\partial^{\alpha-\beta}u\Big\|_{L^2(\R^n)}
\\
   &\leq&  C\sum\limits_{|\alpha|\leq k}\Big\|\partial^{\alpha}u\Big\|_{L^2(\R^n)}+C\sum
   \limits_{|\alpha|\leq k-1}\Big\|\partial^{\alpha}u\Big\|_{L^2(\R^n)}
   \\
   &\leq &
   C\|u\|_{H^k(\R^n)}\;.
\end{eqnarray*}
From $(\ref{eq4.20})$ we deduce that; $\forall\alpha$, $|\alpha|\leq
k$,
$$
e^{-\delta t}\|(\partial^{\alpha}u)(t,\cdot)\|^2_{L^2(\R^n)}\leq
\|(\partial^{\alpha} u)(0,\cdot)\|^2_{L^2(\R^n)}+
2\int_0^te^{-\delta
s}\big\|(\partial^{\alpha}f(s,\cdot)\big\|^2_{L^2(\R^n)}ds+2C\int_0^te^{-\delta
s}\|u(s,\cdot)\|^2_{H^k(\R^n)}ds \;.
$$
Let us take the sum over  $|\alpha|\leq k$; we obtain :
$$
e^{-\delta t}\|u(t,\cdot)\|^2_{H^k(\R^n)}\leq
\|u(0,\cdot)\|^2_{H^k(\R^n)}+ 2\int_0^te^{-\delta
s}\big\|f(s,\cdot)\big\|^2_{H^k(\R^n)}ds+2C\int_0^te^{-\delta
s}\|u(s,\cdot)\|^2_{H^k(\R^n)}ds \;;
$$
 and using  Gronwall's Lemma, we have :
 $$
 e^{-\delta
t}\|u(t,\cdot)\|^2_{H^k(\R^n)}\leq
\Big(\|u(0,\cdot)\|^2_{H^k(\R^n)}+2\int_0^te^{-\delta
s}\big\|f(s,\cdot)\big\|^2_{H^k(\R^n)}ds\Big)e^{2Ct}ds \;.
$$
Note that $\delta$ and $C$ are two positive constants which depend
only on the bounds of the $C^k$-norms of $a$ and $b$. Now setting
$\delta_0=\delta+2C$ and $C_0=2e^{2CT}$ gives
$$
e^{-\delta_0 t}\|u(t,\cdot)\|^2_{H^k(\R^n)}\leq
\|u_0\|^2_{H^k(\R^n)}+C_0\int_0^te^{-\delta_0s}\|f(s,\cdot)\big\|^2_{H^k(\R^n)}ds
\;;
$$
which is the desired inequality. \qed

The weighted version of Proposition \ref{ta} can be stated as
follows:
\begin{cor} \label{c1}Let $d$ be a non negative real number and $k$ an integer.
Under the hypotheses of Proposition \ref{ta} assume further that the
collection of functions $a=\left(a_1, \ldots, a_n\right)$ does not
depend on the space variables ($a\equiv a(t)$) and that $b\equiv 0$.
If $u$ is a $C^1-$solution of the Cauchy problem
(\ref{eq4.1})-(\ref{eq4.2}), then for every $T>0$, if  $f\in
C\big([0,T];H^k_d(\R^n)\big)$,   $u_0\in
 H^k_d(\R^n)$
 then for
all $t\in[0,T], \; u(t, \cdot)\in H^k_d(\R^n)$ and
\begin{equation}\label{eq4.22}
e^{-\delta_1t}\|u(t,\cdot)\|^2_{H^k_d(\R^n)}\leq
\|u_0\|^2_{H^k_d(\R^n)}+C_1\int_0^te^{-\delta_1s}\|f(s,\cdot)\big\|^2_{H^k_d(\R^n)}ds
\end{equation}
where $\delta_1$ and  $C_1$ are positive constants which depend only
on  $\kappa$.
\end{cor}
\pr Let $\beta\in\N^n$ be given such that $|\beta|\leq k$. If we
differentiate Equation (\ref{eq4.1}) with $\partial^{\beta}$ (recall
$a\equiv a(t)$ and
 $b\equiv 0$) and then multiply the differentiated equation by the
 weight $(1+|x|)^{d+|\beta|}$, then we obtain
$$
  \big((1+|x|)^{d+|\beta|}\partial^{\beta}u\big)_t+\sum\limits_{i=1}^na_i
  \big((1+|x|)^{d+|\beta|}\partial^{\beta}u\big)_{x_i}
  =
   (1+|x|)^{d+|\beta|}\partial^{\beta}f +\sum\limits_{i=1}^n(d+|\beta|)\frac{x_i}{|x|}(1+|x|)^{d+|\beta|-1}
   a_i\partial^{\beta}u
  \;.
$$
But the l.h.s is defined by $L$; so we have:
\begin{equation}\label{eq 32''}
L\big[(1+|x|)^{d+|\beta|}\partial^{\beta}u\big]=(1+|x|)^{d+|\beta|}\partial^{\beta}f
 +\sum\limits_{i=1}^n(d+|\beta|)\frac{x_ia_i}{|x|(1+|x|)}(1+|x|)^{d+|\beta|}\partial^{\beta}u
 \;.
\end{equation}
Since the function  $a$ is bounded and $\frac{|x_i|}{|x|(1+|x|)}\leq
1$, we have:
\begin{equation}\label{eq33}
       \bigg\|\sum\limits_{i=1}^n(d+|\beta|)\frac{x_ia_i}{|x|(1+|x|)}
       (1+|x|)^{d+|\beta|}\partial^{\beta}u\bigg\|_{L^2(\R^n)}\leq
       C\|(1+|x|)^{d+|\beta|}\partial^{\beta}u\|_{L^2(\R^n)}\leq
C\|u\|_{H^k_d(\R^n)}
    \;.
\end{equation}
Now, applying Proposition \ref{p1} to  Equation (\ref{eq 32''})
shows that there exists a constant $\delta>0$ such that:
\begin{eqnarray}\label{eq 34}
\lefteqn{e^{-\delta
t}\|(1+|x|)^{d+|\beta|}\partial^{\beta}u(t,\cdot)\|^2_{L^2(\R^n)}\leq}\\
&&\|(1+|x|)^{d+|\beta|}\partial^{\beta}u_0\|^2_{L^2(\R^n)}+2\int_0^te^{-\delta
s}\|f(s,\cdot)\big\|^2_{H^k_d(\R^n)}ds+2C\int_0^te^{-\delta
s}\|u(s,\cdot)\big\|^2_{H^k_d(\R^n)}ds \;;\nonumber
\end{eqnarray}
where we have used (\ref{eq33}). Summing over $|\beta|\leq k$, we
obtain:
\begin{equation}\label{eq 35}
e^{-\delta t}\|u(t,\cdot)\|^2_{H^k_d(\R^n)}\leq
\|u_0\|^2_{H^k_d(\R^n)}+2\int_0^te^{-\delta
s}\|f(s,\cdot)\big\|^2_{H^k_d(\R^n)}ds+2C\int_0^te^{-\delta
s}\|u(s,\cdot)\big\|^2_{H^k_d(\R^n)}ds \;;
\end{equation}
from where using Gronwall's Lemma, we get :
$$e^{-\delta t}\|u(t,\cdot)\|^2_{H^k_d(\R^n)}\leq
\Big(\|u_0\|^2_{H^k_d(\R^n)}+2\int_0^te^{-\delta
s}\|f(s,\cdot)\big\|^2_{H^k_d(\R^n)}ds\Big)e^{2Ct} \;.
$$
Note that $\delta$ and $C$ only depends on the $L^{\infty}$ norm of
$a$ and therefore there exists two positive constants $\delta_1$ and
$C_1$  which depends only on  $\kappa$ such that
$$%
e^{-\delta_1t}\|u(t,\cdot)\|^2_{H^k_d(\R^n)}\leq
\|u_0\|^2_{H^k_d(\R^n)}+C_1\int_0^te^{-\delta_1s}\|f(s,\cdot)\big\|^2_{H^k_d(\R^n)}ds
\;.
$$
\qed
\section{Proof of Lemma \ref{ldd}}

\begin{lemma}
  The space $\mathcal{C}^{\infty}_c(\R^3)$ of Compactly supported smooth functions defined on
  $\R^3$ is dense in the space $H^3_d(\R^3)$.
\end{lemma}
\pr The proof follows closely that of the density of
$\mathcal{C}^{\infty}_c(\R^3)$ in the usual Sobolev space
$H^k(\R^3)$ and will be done in two steps.
\subsubsection*{First step:}   Set $\mathcal{T}=\{f\in
H^3_d(\R^3);\; supp ($f$) \; compact  \}$ and let us show that
 $\mathcal{T}$ is dense in $H^3_d(\R^3)$.
Let $\varphi\in  \mathcal{C}^{\infty}_c(\R^3)$ such that
$
 \varphi (x)=1 $ for
$|x|\leq 1;\; 0\leq\varphi \leq 1$ and $supp(\varphi)\subset
B_{\R^3}(0, 2)$. For every integer $j\geq 1$ set $ \varphi_j(x)=
\varphi\left(\frac{x}{j}\right)\;. $ Then,
$$
\varphi_j\in \mathcal{C}^{\infty}_c(\R^3), \;0\leq\varphi_j \leq
1,\; supp \,\varphi_j\subset B(0,2j) \;.
$$
Moreover, $\,\forall \alpha\in\N^3, \; D^{\alpha}\varphi_j$ is
uniformly bounded; more precisely:
$$
\forall j\in\N^*,\;|D^{\alpha}\varphi_j|\leq
\frac{C_{\alpha}}{j^{|\alpha|}}\leq C_{\alpha}\;\mbox{where}\;
C_{\alpha}=\sup\limits_{x\in\R^3}|D^{\alpha}\varphi(x)|\,.
$$
Since $\varphi_j(x)=1$ for $|x|\leq j$, the sequence $(\varphi_j)_j$
converges pointwise to $1$ as $j$ goes to $\infty$.
Now let $f\in H^3_d(\R^3)$ and set $f_j = \varphi_jf$. Since
$\varphi_j\in\mathcal{C}^{\infty}_c(\R^3),\; f_j\in \mathcal{T}$.
Let us prove that  $(f_j)$ converge towards $f$ in
$H^3_d(\R^3)$.\\
We have
\be\label{e1}
\|f-f_j\|_{H^3_d(\R^3)}=\|f-\varphi_jf\|_{H^3_d(\R^3)}=\Bigg(\sum\limits_{|\alpha|\leq
3}\int_{\R^3}(1+|x|)^{2d+2|\alpha|}\big|D^{\alpha}\big(f(1-\varphi_j)\big)\big|^2dx\Bigg)^{\frac{1}{2}}
\;. \ee
From Leibnitz's formula and the convexity of the function $t\mapsto
t^2$,
$$
\big|D^{\alpha}\big(f(1-\varphi_j)\big)\big|\leq\sum\limits_{|\beta|\leq|\alpha|}
\complement^{\beta}_{\alpha} |D^{\alpha-\beta}f|\,|
D^{\beta}(1-\varphi_j)| \leq
\Bigg(\sum\limits_{|\beta|\leq|\alpha|}\big(\complement^{\beta}_{\alpha}\big)^2\Bigg)^{\frac{1}{2}}
\Bigg(\sum\limits_{|\beta|\leq|\alpha|}|D^{\alpha-\beta}f|^2|D^{\beta}(1-\varphi_j)|^2\Bigg)^{\frac{1}{2}}
\;;
$$
thus,
\be\label{e2}
(1+|x|)^{2d+2|\alpha|}\big|D^{\alpha}\big(f(1-\varphi_j)\big)\big|^2\leq
C_{\alpha}\sum\limits_{|\beta|\leq|\alpha|}(1+|x|)^{2d+2|\alpha|}|D^{\alpha-\beta}f|^2|D^{\beta}(1-\varphi_j)|^2
\;. \ee
Note that
\bea
\sum\limits_{|\beta|\leq|\alpha|}(1+|x|)^{2d+2|\alpha|}|D^{\alpha-\beta}f|^2|D^{\beta}(1-\varphi_j)|^2
&=&(1+|x|)^{2d+2|\alpha|}|D^{\alpha}f|^2|1-\varphi_j|^2\\
&&\;\;\;\;+\sum\limits_{1\leq|\beta|\leq|\alpha|}(1+|x|)^{2d+2|\alpha|}|D^{\alpha-\beta}f|^2|D^{\beta}(1-\varphi_j)|^2
\;. \eea
\textbf{i)}  $\varphi_j\longrightarrow 1$ and
$|1-\varphi_j|\leq1+|\varphi_j|\leq 2$, thus,
$$
\left\{
  \begin{array}{ll}
    (1+|x|)^{2d+2|\alpha|}|D^{\alpha}f|^2|1-\varphi_j|^2\longrightarrow 0,\;\mbox{as}\;j\longrightarrow +\infty \;\mbox{pointwise}
    \\\mbox{and} \\
    (1+|x|)^{2d+2|\alpha|}|D^{\alpha}f|^2|1-\varphi_j|^2\leq 4(1+|x|)^{2d+2|\alpha|}|D^{\alpha}f|^2
  \end{array}
\right. ;
$$
but $f\in H^3_d(\R^3)$, thus $(1+|x|)^{2d+2|\alpha|}|D^{\alpha}f|^2$
is integrable and we can thus apply the Lebesgue dominated
convergence theorem to obtain
$$
\int_{\R^3}(1+|x|)^{2d+2|\alpha|}|D^{\alpha}f|^2|1-\varphi_j|^2dx\longrightarrow
0 \;\mbox{as}\;j\longrightarrow +\infty.
$$
\textbf{ii)} Recall that
$|D^{\beta}\varphi_j|\leq\frac{C_{\beta}}{j^{|\beta|}}, \,\forall
|\beta|\geq 1$ and  $ supp \varphi_j\subset B(0,2j)$,  thus,
$$
(1+|x|)^{2|\beta|}|D^{\beta}\varphi_j|^2\leq
C^2_{\beta}\frac{(1+2j)^{2|\beta|}}{j^{2|\beta|}}\leq
C_{\beta}\frac{(1+j^{2|\beta|})}{j^{2|\beta|}}\leq 2C_{\beta}
$$
and then,
$$
\left\{
  \begin{array}{ll}
    (1+|x|)^{2d+2|\alpha|}|D^{\alpha-\beta}f|^2|D^{\beta}\varphi_j|^2
    \leq(1+|x|)^{2d+2|\alpha|}|D^{\alpha-\beta}f|^2\frac{C_{\beta}}{j}
 \longrightarrow 0,\;\mbox{as}\;j\longrightarrow +\infty \\\mbox{and} \\
    (1+|x|)^{2d+2|\alpha|}|D^{\alpha-\beta}f|^2|D^{\beta}\varphi_j|^2\leq
2C_{\beta}(1+|x|)^{2d+2|\alpha-\beta|}|D^{\alpha-\beta}f|^2\,.
  \end{array}
\right.
$$
Again, $f\in H^3_d(\R^3)$, and then
$(1+|x|)^{2d+2|\alpha-\beta|}|D^{\alpha-\beta}f|^2|$ is integrable.
We use once more the dominated convergence theorem and obtain that
$\forall\,\beta$ such that $1\leq|\beta|\leq|\alpha|$,
$$
\int_{\R^3}(1+|x|)^{2d+2|\alpha|}|D^{\alpha-\beta}f|^2|D^{\beta}\varphi_j|^2dx\longrightarrow
0, \;\mbox{as}\;j\longrightarrow +\infty \;.
$$
and thus
$$
\sum_{1\leq|\beta|\leq|\alpha|}\int_{\R^3}(1+|x|)^{2d+2|\alpha|}|D^{\alpha-\beta}f|^2|D^{\beta}\varphi_j|^2dx\longrightarrow
0 ,\;\mbox{as}\;j\longrightarrow +\infty \;.
$$
Finally we obtain from (\ref{e2}) that
$$
\forall\alpha,\;|\alpha|\leq
3,\;\int_{\R^3}(1+|x|)^{2d+2|\alpha|}\big|D^{\alpha}\big(f(1-\varphi_j)\big)\big|^2dx\longrightarrow
0 \;\mbox{as}\;j\longrightarrow +\infty
$$
and conclude from (\ref{e1}), that $f_j\longrightarrow f$ in the
topology of $H^3_d(\R^3)$. Thus, $\mathcal{T}$ is dense in
$H^3_d(\R^3)$.
\subsubsection*{Second step:}   We prove that
$\mathcal{C}^{\infty}_c(\R^3)$ is dense in $\mathcal{T}$ endowed
with the topology of $H^3_d(\R^3)$. Let $\theta$ be a standard
mollifier, which means that $\theta$ is a positive  $C^{\infty}$
function in $\R^3$ supported in the unit ball and such that
$\int_{\R^3}\theta(x)dx=1$. For every $j\in \N^*$ , we define
$\theta_j$ by $$ \theta_j(x)= j^3\theta(jx) \;.
$$
Then, the sequence $(\theta_j)_{j\in\N^{\ast}}$ has the following
properties%
$$
     \theta_j\geq 0;\;\;\theta_j\in\mathcal{C}^{\infty}_c(\R^3);\;\;
    \int_{\R^3}\theta_j(x)dx=1;\;\;
    supp\,\theta_j\subset\bar{B}(0,\frac1j),\;j\in\N^*
    \;.
$$
Now let $f\in \mathcal{T}$, and let $R$ be a positive real number
such that $supp(f)\subset B(0,R)$. Consider the following
convolution product $f_j=f\star\theta_j$, then $ supp(f_j) $ is
compact since $f$ and $\theta_j$ are compactly supported, more
precisely, $\theta_j\in \mathcal{C}^{\infty}_c(\R^3).$ We shall
prove that $(f_j)$ converges towards  $f$ in $H^3_d(\R^3)$. For this
purpose, it will be sufficient to prove that
$$
\forall\alpha,\; |\alpha|\leq
3,\;(1+|x|)^{d+|\alpha|}D^{\alpha}f_j\longrightarrow
(1+|x|)^{d+|\alpha|}D^{\alpha}f\;\mbox{in}\;L^2(\R^3) \;.
$$
Let  $\alpha,\; |\alpha|\leq 3$, we have:
$D^{\alpha}f_j=D^{\alpha}f\star\theta_j\;\mbox{and}\;D^{\alpha}f-D^{\alpha}f_j=D^{\alpha}f-D^{\alpha}f\star\theta_j$.
Now
$$
D^{\alpha}f\star\theta_j(x)=\int_{\R^3}D^{\alpha}f(y)\theta_j(x-y)dy
$$
and
$$
D^{\alpha}f(x)=D^{\alpha}f(x)\int_{\R^3}\theta_j(x-y)dy=\int_{\R^3}D^{\alpha}f(x)\theta_j(x-y)dy
$$
thus
\be\label{e3}
(1+|x|)^{2d+2|\alpha|}\left|D^{\alpha}f(x)-D^{\alpha}f_j(x)\right|^2\leq
(1+|x|)^{2d+2|\alpha|}\left(\int_{\R^3}
\big|D^{\alpha}f(x)-D^{\alpha}f(y)\big|\theta_j(x-y)dy\right)^2 \;.
\ee
 Applying Holder's inequality gives the following estimates
$$ (1+|x|)^{2d+2|\alpha|}\big|D^{\alpha}f(x)- D^{\alpha}f_j(x)\big|^2
\leq(1+|x|)^{2d+2|\alpha|}\int_{\R^3}\big|D^{\alpha}f(x)-D^{\alpha}f(y)\big|^2\theta_j(x-y)dy\,.
$$
which after integration on $\R^3$, gives:
\bel{e4}
\int_{\R^3}(1+|x|)^{2d+2|\alpha|}\big|D^{\alpha}f(x)-D^{\alpha}f_j(x)\big|^2dx
\leq\int_{\R^3}\int_{\R^3}(1+|x|)^{2d+2|\alpha|}\big|D^{\alpha}f(x)-D^{\alpha}f(y)\big|^2\theta_j(x-y)dxdy\,.
\ee
Consider on $\R^{2n}$ the change of variables $u= x - y , \;v = y$,
we have:
$$
\int_{\R^3}(1+|x|)^{2d+2|\alpha|}\big|D^{\alpha}f(x)-D^{\alpha}f_j(x)\big|^2dx
\leq\int_{|u|\leq\frac{1}{j}}\int_{\R^3}(1+|v+u|)^{2d+2|\alpha|}\big|D^{\alpha}f(v+u)-D^{\alpha}f(v)\big|^2\theta_j(u)dudv\;,
$$
which can be written as
\be\label{e5}
\int_{\R^3}(1+|x|)^{2d+2|\alpha|}\big|D^{\alpha}f(x)-D^{\alpha}f_j(x)\big|^2dx
\leq\int_{|u|\leq\frac{1}{j}}\theta_j(u)du\int_{|v|<R+1}(1+|u+v|)^{2d+2|\alpha|}\big|D^{\alpha}f(u+v)-D^{\alpha}f(v)\big|^2dv\,.
\ee
 then,
for $|u|\leq\frac{1}{j}<1$, we have:
\be\label{e6}
\int_{\R^3}(1+|u+v|)^{2d+2|\alpha|}\big|D^{\alpha}f(u+v)-D^{\alpha}f(v)\big|^2dv\leq
C_{\alpha}\int_{\R^3}\big|D^{\alpha}f(u+v)-D^{\alpha}f(v)\big|^2dv
\ee
where $C_{\alpha}=(3+R|)^{2d+2|\alpha|}\,. $ Recall that $f\in
\mathcal{T}\subset H^3_d(\R^3)\subset H^3(\R^3)$ thus
$D^{\alpha}f\in
 L^2(\R^3)$ and by the  continuity of the $L^2$-norm, one has:
$$
\int_{\R^3}\big|D^{\alpha}f(y+z)-D^{\alpha}f(y)\big|^2dy\longrightarrow
0\;\, \mbox{as}\;\, z\longrightarrow 0
$$
and it follows from (\ref{e6}) that:
$$
\int_{\R^3}(1+|y+z|)^{2d+2|\alpha|}\big|D^{\alpha}f(y+z)-D^{\alpha}f(y)\big|^2dy\longrightarrow
0\;\, \mbox{as}\;\, z\longrightarrow 0 \;.
$$
Now let $ \varepsilon>0$, there exists $\exists\delta>0$ such that
$$
\forall z\in\R^3,
\;|z|\leq\delta\Longrightarrow\int_{\R^3}(1+|y+z|)^{2d+2|\alpha|}\big|D^{\alpha}f(y+z)-D^{\alpha}f(y)\big|^2dy<\varepsilon
\;.
$$
Since  $\lim\limits_{j\rightarrow +\infty}\frac{1}{j}=0$, there
exists $j_0>0$ such that $\forall j\in\N^*$,
$j>j_0\Longrightarrow\frac{1}{j}<\delta$ thus, $j>j_0$ implies
 $|z|<\frac{1}{j}<\delta$ and then
$$
\int_{\R^3}(1+|y+z|)^{2d+2|\alpha|}\big|D^{\alpha}f(y+z)-D^{\alpha}f(y)\big|^2dy<\varepsilon
\;;
$$
from where we obtain
$$
j>j_0\Longrightarrow\int_{\R^3}(1+|x|)^{2d+2|\alpha|}\big|D^{\alpha}f(x)-D^{\alpha}f_j(x)\big|^2dx
\leq\varepsilon\int_{|z|\leq\frac{1}{j}}\theta_j(z)dz<\varepsilon\,;
$$
thus,
$$
\forall\alpha,\; |\alpha|\leq
3,\;(1+|x|)^{d+|\alpha|}D^{\alpha}f_j\longrightarrow
(1+|x|)^{d+|\alpha|}D^{\alpha}f\;\mbox{in}\;L^2(\R^3) \;;
$$
and consequently, $f_j\longrightarrow f$ in $H^3_d(\R^3)$ which
proves that $\mathcal{C}^{\infty}_c(\R^3)$ is dense in $\mathcal{T}$
endowed with the topology of $H^3_d(\R^3)$.
  \qed

\bibliographystyle{amsplain}

\addcontentsline{toc}{1}{\bigskip \noindent  {\bf References
\hfill}}

\bibliography{Biblio-NK}

\def\cprime{$'$} \def\cprime{$'$} \def\cprime{$'$} \def\cprime{$'$}
\providecommand{\bysame}{\leavevmode\hbox
to3em{\hrulefill}\thinspace}
\providecommand{\MR}{\relax\ifhmode\unskip\space\fi MR }
\providecommand{\MRhref}[2]{%
  \href{http://www.ams.org/mathscinet-getitem?mr=#1}{#2}
} \providecommand{\href}[2]{#2}

\end{document}